\documentclass[12pt]{amsart}
\usepackage{amsmath, amssymb}
\newfont {\cyr} {wncyr10}
\renewcommand{\labelenumi}{{(\roman{enumi})}}
\usepackage{color}
\newcommand{\red}[1]{\,{\color{black} #1}\,}
\newcommand{\blue}[1]{\,{\color{black} #1}\,}
\usepackage{amsfonts}
\usepackage{amsmath}
\usepackage{amssymb}
\usepackage{setspace}
\usepackage{a4}
\usepackage{multicol}

\newtheorem{theorem}{Theorem}[section]
\newtheorem{lemma}[theorem]{Lemma}
\newtheorem{proposition}[theorem]{Proposition}
\newtheorem{corollary}[theorem]{Corollary}
\newtheorem{definition}[theorem]{Definition}

\newcounter{claim}[theorem]

\newcounter{cclaim}[theorem]

\def \signal {\mbox{{{\cyr I}}}}

\def \t {\tau}

\newcommand{\E}{\mathrm{E}}\newcommand{\SU}{\mathrm{SU}}
\newcommand{\F}{\mathrm{F}}\newcommand{\M}{\mathcal{M}}
\newcommand{\G}{\mathrm{G}}

\newcommand{\C}{\mathcal{C}}
\newcommand{\Q}{\mathrm{Q}}

\newcommand{\Aut}{\mathrm{Aut}}

\newcommand{\Out}{\mathrm{Out}}

\newcommand{\Syl}{\mathrm{Syl}}\newcommand{\syl}{\mathrm{Syl}}

\newcommand{\GF}{\mathrm{GF}}
\newcommand{\GL}{\mathrm{GL}}
\newcommand{\Sp}{\mathrm{Sp}}
\newcommand{\SL}{\mathrm{SL}}

\newcommand{\PSL}{\mathrm{PSL}}\newcommand{\PSp}{\mathrm{PSp}}
\newcommand{\Sym}{\mathrm{Sym}}
\newcommand{\Alt}{\mathrm{Alt}}
\newcommand{\Dih}{\mathrm{Dih}}
\newcommand{\SDih}{\mathrm{SDih}}

\newcommand{\U}{\mathrm{U}}

\def \GU {\mbox {\rm GU}}

\def \syl {\hbox {\rm Syl}}\def \Syl {\hbox {\rm Syl}}

\def \ov {\overline}

\def \wt {\widetilde}

\def \Aut{ \mathrm {Aut}}

\def \Out{\mbox {\rm Out}}

\def \M{\mbox {\rm M}}

\def \Co {\mbox {\rm Co}}

\def \McL{\mbox {\rm McL}}

\def \PSU {\mbox {\rm PSU}}
\def \GSp {\mbox {\rm GSp}}\def \GO {\mbox {\rm GO}}

\begin{document}
\renewcommand{\labelenumi}{(\roman{enumi})}

\title  {An identification theorem for ${\mathbf \PSU_6(2)}$ and its automorphism groups}
 \author{Chris Parker}
  \author{Gernot Stroth}

\address{Chris Parker\\
School of Mathematics\\
University of Birmingham\\
Edgbaston\\
Birmingham B15 2TT\\
United Kingdom} \email{c.w.parker@bham.ac.uk}

\address{Gernot Stroth\\
Institut f\"ur Mathematik\\ Universit\"at Halle - Wittenberg\\
Theordor Lieser Str. 5\\ 06099 Halle\\ Germany}
\email{gernot.stroth@mathematik.uni-halle.de}

\email {}

\date{\today}

\maketitle \pagestyle{myheadings}

\markright{{\sc $\PSU_6(2)$ and its automorphism groups }} \markleft{{\sc Chris Parker and Gernot Stroth}}

\begin{abstract} We identify the groups $\PSU_6(2)$, $\PSU_6(2){:}2$, $\PSU_6(2){:}3$ and $\Aut(\PSU_6(2))$ from the  structure of the centralizer of an element of order $3$.
 \end{abstract}

\section{Introduction}

\blue{
The aim of this article is to provide a $3$-local identification of the groups $\PSU_6(2){:}3$, $\PSU_6(2){:}2$ and  $\PSU_6(2){:} \Sym(3)$  as such characterizations  are needed in the ongoing work to classify groups with a large $p$-subgroup. For a prime $p$, a \emph{$p$-local} subgroup of $G$ is by definition the normalizer in $G$ of a non-trivial $p$-subgroup of $G$. We say that a  $p$-subgroup $Q$ of a group $G$ is \emph{large}
 provided
\medskip
\begin{enumerate}
\item[(L1)]\label{1} $F^*(N_G(Q))=Q$; and
\item[(L2)]\label{2} if $U$ is a non-trivial subgroup of $ Z(Q)$, then $N_G(U)\le N_G(Q)$.
\end{enumerate}
\medskip

An interesting  observation is that most of the groups of  Lie type in characteristic $p$ contain a  \emph{large}
$p$-subgroup. In fact the only Lie type groups in characteristic $p$ and rank at least $2$ which do not contain
such a subgroup are $\PSp_{2n}(2^a)$, $\F_4(2^a)$ and $\G_2(3^a)$. It is not difficult to show that groups $G$
which contain a large $p$-subgroup are of parabolic characteristic $p$, which means that all $p$-local overgroups
$N$ of a Sylow $p$-subgroup $S$  satisfy $F^*(N)=O_p(N)$ (\cite[Lemma 2.1]{PS1}).  The work initiated in
\cite{MSS} starts the determination of the $p$-local overgroups of $S$ which are not contained in $N_G(Q)$. This
is the first mile of a long road to showing that typically a group with a large $p$-subgroup is a group of Lie
type defined in characteristic $p$ and of rank at least $2$. The basic crude idea is to gather information about
the $p$-local subgroups of $G$ containing a fixed Sylow $p$-subgroup  so that the subgroup generated by them can
be identified with a  group of Lie type via its action on the chamber complex coming from these subgroups (which
will in fact be the maximal parabolic subgroups).  However,
 one is  sometimes confronted with the following situation.
Some (but perhaps not all) of the  $p$-local subgroups of $G$ containing a given Sylow $p$-subgroup $S$ of $G$
generate a subgroup   $H$ and $F^*(H)$  is known to be isomorphic to a Lie type group in characteristic $p$. The
expectation (or rather hope) is that $G=H$. In the case that $H$ is a proper subgroup of $G$,   one usually tries
to prove that $H$ contains \emph{all} the  $p$-local subgroups of $G$ which contain $S$ and then in a next step
to prove that $H$ is strongly $p$-embedded in $G$ at which stage \cite{PSStrong} is applicable and delivers $G =
H$. The last two steps are reasonably well understood, at least for groups with mild extra assumptions. However
it might be that the first step cannot be made.  Typically this will occur only when $N_G(Q)$ is not contained in
$H$.
 If $N_H(Q)$  is non-soluble and $p$ is odd,  A. Seidel in his PhD thesis
\cite{Seidel} has shown that this cannot occur. In  \cite{PS1} the authors use the identification theorem
presented in this paper together with further identifications \cite{McL,ParkerRowley2,  PS,F42,PSS}  to handle
the more delicate analysis when $p=3$ and $N_H(Q)$ is soluble. Far from these configurations not arising, the
rule of thumb here is that if it might happen then it does. The possibilities for the group $F^\ast(H)$ are
easily shown to be $\PSL_3(3^a)$ (which we do not consider), $\PSp_4(3)$, $\PSU_4(3)$, $\PSL_4(3)$,
$\Omega_7(3)$, $\mathrm P\Omega^+_8(3)$ and $\G_2(3^a)$. The main theorem in \cite{PS1} says that if $N_G(Q) \not
\le H$, then  $F^*(G)$ is one of  $\U_6(2)$, $\F_4(2)$, ${}^2\E_6(2)$, $\McL$,  $\Co_2$, $\M(22)$,  $\M(23)$ or
$\F_2$.

Thus the present article is required for the identification of $\U_6(2)$ and its automorphism groups in the  case when $F^*(H) \cong \U_4(3)$. Furthermore, this article plays a pivotal role in the $3$-local identifications
of $\M(22)$ and ${}^2\E_6(2)$. Indeed the centralizers of involutions in both  $\M(22)$ and ${}^2\E_6(2)$ and their automorphism groups feature $\U_6(2)$ and its automorphism groups prominently.  These identifications  in addition are required for    our  work on the sporadic simple groups $\M(23)$ and $\F_2$.}

In earlier work \cite{Parker1} the first author  proved  the
following result: let $G$ be a finite group, $S$ be a Sylow
3-subgroup of $G$ and $Z = Z(S)$. Assume that $N_G(Z)$ is similar to
a 3-normalizer in $\PSU_6(2)$. Then either $Z$ is weakly closed in
$S$ or $G \cong \PSU_6(2)$.  However,  for our intended applications of such
results as outlined above,  we also  need to  identify the groups
$\PSU_6(2){:}3$, $\PSU_6(2){:}2$ and  $\PSU_6(2){:} \Sym(3)$ from
their $3$-local data  (here and throughout this work  we use
the  Atlas \cite{Atlas} notation for group extensions). The addition of these automorphisms cause numerous  difficulties.

\red{\begin{definition} We say that  $X$ is similar to a $3$-centralizer in a group of type  $\PSU_6(2)$ or $\F_4(2)$
provided the following conditions hold.
\begin{enumerate}
\item $Q=F^*(X)$ is extraspecial of order $3^5$ and $Z(F^*(X)) =Z(X)$; and
\item $X/Q$ contains a normal subgroup isomorphic to $ \Q_8\times \Q_8$.
\end{enumerate}
\end{definition}}

A precise description of the possibilities for the group $X/Q$ will be determined in Section~3.
Our theorem is as follows.

\begin{theorem}\label{MainTheorem} Suppose that $G$ is a group, $Z \le G$ has order $3$ and set $H= C_G(Z)$. If $H$ is similar to
a $3$-centralizer in a  group of type $\PSU_6(2)$  of $\F_4(2)$ and $Z$ is  weakly
closed in $F^*(H)$ but not in $H$, then $G \cong \PSU_6(2)$,
$\PSU_6(2){:}2$, $\PSU_6(2){:}3$ or $\PSU_6(2){:}\Sym(3)$.
\end{theorem}

In the case that $Z$ is weakly closed in $H$, then $G$ could be a nilpotent group
extended by a group similar to a $3$-centralizer
 of type $\PSU_6(2)$ of $\F_4(2)$. Thus the hypothesis that $Z$ is not weakly
closed in  $H$ is necessary to have an identification theorem. On
the other hand, the hypothesis that $Z$ is weakly closed in $F^*(H)$
is there to prevent further examples related to $\F_4(2)$ arising.
The methods that we use here are also be applicable
to this type of configuration, however the investigation of such a
possibility would take a rather different road at the very outset of
our proof and so the analysis of this possibility is not
included here and is the subject of \cite{F42}. Combining the work of both papers we obtain

\red{
\begin{theorem}\label{combT} Suppose that $G$ is a group, $Z \le G$ has order $3$ and set $H = C_G(Z)$.
If $H$ is similar to a $3$-centralizer of a  group of type $\PSU_6(2)$ or $\F_4(2)$ and $Z$ is not weakly closed in
a Sylow $3$-subgroup of $G$  with respect to $G$, then  either $F^*(G) \cong \F_4(2)$ or  $F^*(G) \cong \PSU_6(2)$.
\end{theorem}
}
We now describe the layout of the paper and highlight a number of
interesting features of the article.    We begin in Section~2 with
preliminary lemmas and background material. Noteworthy results in
this section are Lemma~\ref{H+}  where we embellish  the statement
of Hayden's Theorem \cite{Hayden} to give the structure of the
normal subgroup of index $3$ and Lemma~\ref{notsimple} where we use
transfer theorems to show that a group with a certain specified
$2$-local subgroup  has a subgroup of index $2$.  The relevance of
such results to our proof is apparent as a look at the list of
groups in the conclusion of our theorem shows. Let $G$, $H$  and $Z$  be as in the
statement  of Theorem~\ref{MainTheorem} and let $S\in \Syl_3(M)$ where $M=N_G(Z)$ contains $H$ at index at most $2$. In
Section~3, we tease out the structure of $M$ and establish much of
the notation that is used throughout the proof of
Theorem~\ref{MainTheorem}.

 In Section~4, we determine the structure
the normalizer of a further $3$-subgroup which we call $J$ and turns
out to be the Thompson subgroup of $S$. The fact that $N_G(J)$ is
not contained in $M$ is a consequence of the hypothesis that  $Z$ is
not weakly closed in $M$. We find in Lemma~\ref{S6}  that $N_G(J)/J
\cong 2\times \Sym(6)$ or $\Sym(6)$. With this information, after using a transfer theorem, we are  able
to apply \cite{Parker1} and do so  in  Theorem~\ref{NMSagain} to get that $G
\cong \PSU_6(2)$ or $\PSU_6(2){:}3$ if $N_M(S)/S \cong \Dih(8)$.
Thus from this stage on we assume  that $N_M(S)/S \cong 2 \times
\Dih(8)$ and $N_M(J)/J \cong 2 \times \Sym(6)$. With this
assumption,  our target groups  all have a subgroup of index $2$.
Our plan is to determine the structure of a $2$-central involution
$r$, apply Lemma~\ref{notsimple} and then apply
Theorem~\ref{NMSagain} to the subgroup of index $2$. The involution
we focus on is contained in $M$ and centralizes a subgroup of
$F^*(M)$  isomorphic to $3^{1+2}_+$. But before we can make this
investigation we need to determine the centralizer of another
subgroup (for now we will call it $X$) which has order either $3$ or
$9$. It turns out we may apply the theorems of Hayden \cite{Hayden}
and Prince \cite{prince1} to get  $E(C_G(X)) \cong \SU_4(2)$.  At this juncture, given the $3$-local
information that we have gathered, we can  construct an
extraspecial $2$-subgroup  $\Sigma$ of order $2^9$ in  $K=C_G(r)$.
In Theorem~\ref{normalizer} we show that $N_K(\Sigma)/\Sigma \cong
\Aut(\SU_4(2))$, $(\SU_4(2) \times 3){:}2$ or $\Sp_6(2)$.
 In our target groups the possibility $\Sp_6(2)$ does not
arise and we will say more about this shortly.

In Section~6
we show that $\Sigma$ is strongly closed in $N_K(\Sigma)$ with
respect to $K$ and then we apply Goldschmidt's Theorem to get that
$K= N_K(\Sigma)$. At this stage we know the centralizer of a
$2$-central involution and so we prove the theorem in Section~6.  We
mention here that when $K/\Sigma \cong \Sp_6(2)$ we apply
\cite{Smith} to obtain  $G \cong \Co_2$ and then  eliminate
this group as it does not satisfy our hypothesis on  the structure
of $M$. One should wonder if the configuration involving $\Sp_6(2)$
could be eliminated at an earlier stage. However, as $\Co_2$
contains $\PSU_6(2){:}2$ as a subgroup of index $2300$, these groups
are intimately related. A $3$-local identification of $\Co_2$ can be
found in \cite{ParkerRowley2}.

Our notation follows that in \cite{Aschbacher}, \cite{Gorenstein} and  \cite{GLS2}. \blue{In particular we use
the definition of signalizers as given in \cite[Definition 23.1]{GLS2} as well as the notation $\signal_G(A,\pi)$
to denote the set of $A$-signalizers in $G$ and $\signal_G^*(A,\pi)$ the maximal  members of $\signal_G(A,\pi)$.}
As mentioned earlier we use Atlas \cite{Atlas} notation for group extensions. We also use \cite{Atlas} as a
convenient source for information about subgroups of almost simple groups. Often this information can be easily
gleaned from well-known properties of classical groups. For odd $p$, the extraspecial groups of exponent $p$ and
order $p^{2n+1}$ are denoted by $p^{1+2n}_+$. The extraspecial $2$-groups of order $2^{2n+1}$ are denoted by
$2^{1+2n}_+$ if the maximal elementary abelian subgroups have order $2^{1+n}$ and otherwise we write
$2^{1+2n}_-$. We hope our notation for specific groups is self-explanatory. In addition, for a subset $X$ of a
group $G$, $X^G$ denotes that set of $G$-conjugates of $X$. If $x, y \in H \le G$, we often write $x\sim _Hy$ to
indicate that $x$ and $y$ are conjugate in $H$. All the groups in this paper are finite groups.

\medskip

\noindent {\bf Acknowledgement.}  The initial draft of this paper was prepared during  a visit
to the Mathematisches Forschungsinstitut Oberwolfach as part of the
Research in Pairs Programme, 30th November--12 November, 2009. The authors are pleased to thank the MFO and its staff for the pleasant and
stimulating research environment that they provided. The first author is also grateful to the DFG for financial support and the mathematics department in Halle for their hospitality.

\section{Preliminaries}

In this section we gather preliminary results for our proof of Theorem~\ref{MainTheorem}.  For a group $G$
with Sylow $p$-subgroup $P$ and $v \in P$,  $v$ is  said to be \emph{extremal} in $P$ if $C_P(v)$ is a Sylow $p$-subgroup
of $C_G(v)$.

\begin{lemma}\label{Thompson} Suppose that $p$ is a prime and  $G$ is a group. Let  $P$  a Sylow $p$-subgroup of $G$ and $Q$ be a proper normal  subgroup of  $P$ such that $P/Q$ is cyclic. Assume there is $u \in P \setminus Q$ such that
\begin{itemize}
\item[(a)] no conjugate of $u^p$ is contained in $P \setminus Q$; and
\item[(b)] any extremal conjugate of $u$ in $P$ is contained in $Q\cup Qu$.
\end{itemize}
Then either $G$ has a normal subgroup $N$ with $G/N$ cyclic and $u
\not\in N$ or there is $g \in G$ such that
\begin{itemize}
\item[(i)] $u^g \in Q$;
\item[(ii)] $u^g$ is extremal in $P$; and
\item[(iii)] $C_P(u)^g \leq C_P(u^g)$.
\end{itemize}
\end{lemma}

\begin{proof} See \cite[Proposition 15.15]{GLS2} or \cite[Corollary 5.3.1]{Yoshida}.
\end{proof}

\begin{lemma}\label{ptransfer} Suppose that $p$ is a prime, $G$ is a group and $P \in \Syl_p(G)$.
\begin{itemize}
\item[(i)] Assume that there is a normal subgroup $Q$ of $P$ such that $P/Q$ is cyclic and that $y \in P \setminus Q$ has order $p$.
 If every extremal conjugate of $y$ in $P$ is contained in $Qy$, then $G$ has a normal subgroup $N$ with $y \not \in N$ and  $G/N$  cyclic.
\item[(ii)] Assume that $P \leq M \leq G$, $y \in P\setminus M'$ has order $p$ and that, if  $x\in G$ with $y^x \in P$  extremal, then there is $g \in M$ such that $y^x = y^g$. Then $y \not\in G^\prime$.
\item[(iii)] Assume that $J=J(P)$ is the Thompson subgroup of $P$. If  $J$ is elementary abelian  and $J \not \le N_G(J)'$, then $J \not \le G'$.
\end{itemize}
\end{lemma}

\begin{proof} (i) This follows from \ref{Thompson}.

(ii) As $M/M^\prime$ is abelian, there is $N \leq M$ such that
$M^\prime \leq N$, $y \not\in N$, $M = NP$ and $P/(P \cap N)$ is
cyclic. Set $Q = P \cap N$. Now for $g \in M$ with $y^g \in P$ we
have that $y^g \in Qy$. Hence by assumption $y^x \in Qy$ for all $x
\in G$ such that $y^x$ is extremal in $P$. Now (ii) follows from
(i).

(iii) Set $M= N_G(J)$ and pick $y \in J \setminus M'$. Assume that
$g \in G$ and $y^g $ is extremal in $P$.
 Then $C_P(y^g)\in \Syl_p(C_G(y^g))$. Since $C_G(y)$ contains $J$, we have $C_P(y^g)$ contains a $G$-conjugate of $J$. Since $J$ is weakly closed in $P$, we have
 $J \le C_P(y^g)$. But then $y^g \in C_P(J) \le J$. Since $M$ controls fusion in $J$, we now have that $y^g = y^m$ for some $m\in M$. Now (iii) follows from (ii).
\end{proof}

\begin{lemma}\label{quadratic form}
Suppose that $F$ is a field, $V$ is an $n$-dimensional vector space
over $F$ and  $G= \GL(V)$. Assume that $q$ is quadratic form  of
Witt index at least $1$ and $\mathcal S$ is the set of singular
1-dimensional subspaces of $V$ with respect to $q$. Then the
stabiliser in $G$ of $\mathcal S$ preserves $q$ up to similarity.
\end{lemma}

\begin{proof} See \cite[Lemma 2.10]{ParkerRowley2}.
\end{proof}

\begin{lemma}\label{GO4} Suppose that $p$ is an odd prime, $X = \GL_4(p)$ and $V$ is the natural $\GF(p)X$-module. Let $A =\langle a, b\rangle\le X$ be elementary abelian of order $p^2$ and assume that $[V,a] = C_V(b)$ and $[V,b]= C_V(a)$ are distinct and of dimension $2$.
Let $v \in V\setminus [V,A]$. Then  $A$ leaves invariant a
non-degenerate quadratic form with respect to which  $v$ is a
singular vector. In particular, $X$ contains exactly two conjugacy
classes of subgroups such as $A$. One is conjugate to a Sylow
$p$-subgroup of $\GO_4^+(p)$ and the other to a Sylow $p$-subgroup
of $\GO_4^-(p)$.
\end{lemma}
\begin{proof} See \cite[Lemma 2.11]{ParkerRowley2}.
\end{proof}

\begin{lemma}\label{H+} Let $H$ be a finite group and let  $d \in H$ be an element of order $3$ such that $X = C_H(d)$ is isomorphic
to the centralizer  of a non-trivial $3$-central element in
$\PSp_4(3)$.  Let $P \in \Syl_3(X)$ and $E$ be the
elementary abelian subgroup of $P$ of order $27$. Assume that $E$
does not normalize any non-trivial $3^\prime$-subgroup of $H$, that
 $d$ is not $H$-conjugate to its inverse and $H$ has a normal
 subgroup of index $3$. Then $H=C_H(d)$.
\end{lemma}

\begin{proof} Notice first of all that $P \in \Syl_3(H)$. Let $H_1$ be a normal subgroup of $H$ of index $3$
and set $E_1= E\cap H$. So $C_{H_1}(d) \cong 3^{1+2}_+{:}\Q_8$ and
$E_1$ has order $9$. Suppose that $x \in E_1\setminus\langle d
\rangle$. We see that all subgroups of order three in $E_1$ different from $\langle d \rangle$ are conjugate in $O_3(C_H(d))$ and so all $x \in E_1 \setminus \langle d \rangle$ are  conjugate to its inverse and $d$ is not,
$d$ is the unique conjugate of $d$ in $E_1$. Furthermore, $d$ is not
conjugate to any element of $E\setminus H'$ and so $d$ is the unique
conjugate of $d$ in $E$.  Since $x$ is not conjugate to $d$, we have
that $E_1=\langle d,x\rangle $ is a Sylow $3$-subgroup of
$C_{H_1}(x)$. As $E_1/\langle x \rangle $ is self-normalizing in
$C_{H_1}(x)/\langle x\rangle$, $C_{H_1}(x)$ has a normal
$3$-complement $T$ by Burnside's Theorem. However $C_{H_1}(x)$ is
normalized by $E$ and so $T= 1$ by hypothesis. It follows that
$C_H(x)= E$ for all $x \in E_1 \setminus \langle d \rangle$.

Let $y \in E\setminus H_1$. Then, as before,  $E_1$ is a Sylow
$3$-subgroup of $C_{H_1}(y)$. Since $d$ is not conjugate to any
non-trivial element of $E_1 \setminus \{d\}$, we have $N_{H}(E_1)\le
X$. So $N_{C_{H_1}(y)}(E_1)= \langle E_1, s\rangle$ where $s$ is an
element of order at most two in $X$. Since $[E_1,s]< E_1$, Gr\"un's
Theorem \cite[Chapter 7, Theorem 4.4]{Gorenstein}  implies that $C_{H_1}(y)$
has a subgroup $L$ of index at least $|E_1:[E_1,s]|$ with Sylow
$3$-subgroup $[E_1,s]$.  Since $L$ is normalized by $E$, we also
have $O_{3'}(L)=1$. Hence, if $s=1$, then  $C_H(y)\le X$ which means
that $C_H(y)=E$. So suppose that $[E_1,s]$ has order $3$.  Then, as
$C_H([E_1,s])=E$, we have $[E_1,s]$
 is self-centralizing in $L$.  Applying the other Feit-Thompson Theorem
\cite{FeitThompson} to $L$ and using $O_{3'}(L)=1$, we now have that
either $L\cong \Sym(3)$ with $L= N_{X\cap H_1}(
[E_1,s])$ or $L \cong \PSL_3(2)$ or $\Alt(5)$. The latter two
cases are eliminated as $L$ is normalized by $E_1$ and the
centralizers of all of the non-trivial elements of $E_1$ are
soluble. Therefore, $C_H(y)= C_X(y) \le X$ for all $y \in E\setminus
E_1$.

Now let $R \in \Syl_2(X)$ and $r\in R$ be an involution. Then
$C_X(r) = R\langle d,y\rangle$ for some $y\in E\setminus E_1$.
Furthermore, as $d$ is the unique conjugate of $d \in \langle d,
y\rangle$, $$N_{C_H(r)}(\langle d,y\rangle)= N_{X}(\langle
d,y,r\rangle)=\langle d,y,r\rangle$$  and so $C_H(r)$ has a normal
$3$-complement $U$ by Burnside's Theorem. Finally $$U= \langle
C_U(w) \mid w \in \langle d,y\rangle^\#\rangle\le X$$ as $C_H(w) \le
X$ for each $w \in \langle d,y \rangle^\#$. It follows that $U=R$.
But then $R \in \Syl_2(H)$ and $r \in Z^*(H)$ by
\cite{BrauerSuzuki}. As $[O_3(X) ,r]= O_3(X)$, we conclude $O_3(X)
\le O_{2'}(H)$ and deduce $H=X$ from  the Frattini Argument. This
completes the proof of the lemma.
\end{proof}

\begin{theorem} [Hayden]\label{Hayden} Let $H$ be a finite group and let $d$ be an element of order $3$ in $H$ such that $X = C_H(d)$   is isomorphic
to the centralizer  of a non-trivial $3$-central element in
$\PSp_4(3)$. Let $P \in \Syl_3(X)$ and $E$ be the
elementary abelian subgroup of $P$ of order $27$. If  $E$ does not
normalize any non-trivial $3^\prime$-subgroup of $H$ and $d$ is not
$H$-conjugate to its inverse, then either $H=X$ or  $H\cong
\PSp_4(3)$.
\end{theorem}

\begin{proof} By \cite{Hayden}  either $H \cong \PSp_4(3) $ or $H$ has a normal subgroup of index $3$. The result now  follows from Lemma~\ref{H+}.
\end{proof}
\begin{theorem}[A. Prince] \label{PrinceThm}
Let $H$ be a finite group and let $d$ be an element of order $3$ in $H$ such that $X = C_H(d)$   is isomorphic
to the centralizer  of a non-trivial $3$-central element in
$\PSp_4(3)$.  Let $P \in
\Syl_3(C_H(d))$ and $E$ be the elementary abelian subgroup of $P$ of
order $27$. If $E$ does not normalize any non-trivial
$3^\prime$-subgroup of $H$ and $d$ is $H$-conjugate to its inverse,
then either
\begin{enumerate}
\item
$|H:C_H(d)| =2$;
\item $H$ is isomorphic to $\Aut(\SU_4(2))$; or
\item $H$ is isomorphic to $\Sp_6(2)$.
\end{enumerate}
\end{theorem}

\begin{proof} See \cite[Theorem 2]{prince1}. \end{proof}

\begin{lemma} \label{cen3psp43}Suppose that $X$ is a group of shape $3^{1+2}_+.\SL_2(3)$,  $O_2(X)=1$
 and
a Sylow $3$-subgroup of $X$ contains an elementary abelian subgroup of order $3^3$. Then $X$ is isomorphic to the
centralizer of a non-trivial $3$-central element in $\PSp_4(3)$.
\end{lemma}
\begin{proof} See \cite[Lemma~6]{Parker1}.\end{proof}

\begin{lemma}\label{Parker} Let $G$ be a finite group and $S$ be a Sylow 3-subgroup of $G$. Set $Z = Z(S)$, $J = J(S)$ and $M = N_G(Z)$. Suppose that $G^* $ is a  normal subgroup of $G$ and set $M^*=M \cap G^*$. Assume that the
following hold: \begin{enumerate} \item $|M^*| = 2^7.3^6$;
 \item
$M^* \ge QR= O_{3,2}(M^*)$, where $Q = O_3(M^*)$ is extraspecial of order
$3^5$ and $R \in \Syl_2(O_{3,2}(M^*))$; \item $O^2(M^*)= (S\cap M^*)R$ has index $2$ in $M^*$; and
\item $Q/Z$ is a $M^*$-chief factor. \end{enumerate} If
$N_{G^*}(J\cap G^*)\not \le M^*$, then $G^* \cong \PSU_6(2)$ and $G$
is a subgroup of $\Aut(\PSU_6(2))$ such that $G/G^*\cong M/M^*$.
\end{lemma}

\begin{proof} Since  $N_{G^*}(J\cap G^*) \not \le M^*$, $Z$ is not weakly closed in $S \cap G^*$.
The conditions imposed on the structure of $M^*$ mean  that $M^*$ is
similar to a $3$-normalizer in $\PSU_6(2)$ \cite[Definition
1]{Parker1}. Hence \cite[Theorem 1]{Parker1} gives the result.
\end{proof}

\begin{lemma}\label{involutionsonexspec}
Suppose that $E$ is an extraspecial $2$-group and $x \in \Aut(E)$ is
an involution. If $C_E(x)\ge [E,x]$, then $[E,x]$ is elementary
abelian.
\end{lemma}

\begin{proof} Let $\langle e \rangle= Z(E)$.  We show that every element of $[E,x]$ has order $2$. Let $f \in [E,x] \setminus \langle e \rangle$. Then $fe$  has the same order as $f$. Thus we may suppose that $f=[h,x]$ for some $h \in E$.  As $x[h,x]=[h,x]x$ by hypothesis, we have
\begin{eqnarray*}f^2&=&[h,x][h,x] = h^{-1}xhx[h,x]=h^{-1}xh[h,x]x\\&=& h^{-1}xhh^{-1}xhxx=1\\
 \end{eqnarray*}as required. This proves the lemma.
\end{proof}

\red{The following lemma  is an easy consequence of the Three Subgroup Lemma.
 \begin{lemma}\label{w=1} Suppose that $p$ is a prime, $P$ is a $p$-group of nilpotency class at most $2$ and that $\alpha\in \Aut(P)$ has order coprime to $p$.
If $\alpha$ centralizes a maximal abelian subgroup of $P$, then $\alpha=1$.
\end{lemma}

\begin{proof}
See \cite[Lemma~{2.3}]{PS1}.
\end{proof}}

For use in Lemma~\ref{notsimple} and Section~6, we collect some facts
about the action of $\Sp_6(2)$ and $\Aut(\SU_4(2))$ on their
irreducible $8$-dimensional module $V$ over $\GF(2)$. Recall that
$\Aut(\SU_4(2))\cong \mathrm O_6^-(2)$ is a subgroup of $\Sp_6(2)$
\cite[page 46]{Atlas}. We will frequently use the fact that as
$\SU_4(2)$-module, $V$ is the natural $4$-dimensional
$\GF(4)\SU_4(2)$-module regarded as a module over $\GF(2)$. We will
often refer to this as the \emph{natural} $\SU_4(2)$-module.

\begin{table}
\begin{tabular}{|c|c|c|c|}
\hline
&$\Aut(\SU_4(2))$&$\Sp_6(2)$&$\dim C_V(u_j)$\\
\hline $u_1$&$2^{1+4}_+.(\Sym(3) \times \Sym(3))$&$2^7.(\Sym(3)
\times \Sym(3))$&6\\ $u_2$&
$2^6. 3$&$2^7. 3$&4\\
$u_3$&$2 \times \Sym(6)$&$2^5.\Sym(6)$&4\\
$u_4$&$2 \times (\Sym(4) \times 2)$&$2^9. 3$&4\\\hline
\end{tabular}
\caption{Involutions in $\Sp_6(2)$ and $\Aut(\SU_4(2))$. The
involutions in the first row are the \emph{unitary transvections.}
The involutions in the last two rows are those which are in
$\Aut(\SU_4(2)) \setminus \SU_4(2)$.} \label{Table1}
\end{table}

\begin{proposition}\label{fact} Let $X \cong \Sp_6(2)$ and $Y \cong
\Aut(\SU_4(2))$. Assume that $V$ is the $8$-dimensional irreducible
module for $X$ (and hence $Y$)
 over $\GF(2)$. Then the following hold:
\begin{enumerate}
\item $X$ and $Y$ both possess exactly four conjugacy classes of involutions. In
Table~\ref{Table1} we list the four classes of involutions and give
structural information about the centralizers in both groups as
can be found  in \cite[pages 26 and 46]{Atlas}.

\item $X$ and $Y$ have orbits of length $135$ and  $120$  on the non-zero
elements of $V$. We call elements of the orbits non-singular and singular vectors respectively.  Suppose that $x$ is singular and $y$ is
non-singular. Then
\begin{eqnarray*}
|C_{Y}({x})| = 2^7\cdot 3 ,&& |C_{X}({x})| = 2^9\cdot 3 \cdot 7.\\
C_{Y}({y}) \cong 3^{1+2}_+.\SDih(16), &&C_{X}({y}) \cong
\G_2(2).\end{eqnarray*}

\item $X$ and $Y$ both have  exactly three conjugacy classes of
elements of order $3$. They are distinguished  by their action on
$V$. They have centralizers of dimension $0$, $2$ and $4$. The
elements with centralizer  of dimension $2$ are $3$-central and
centralize only non-singular vectors in $V^\#$.

\item For $u\in Y$ an involution, $\dim{C_V}(u)$ is given in column $4$ of
Table~\ref{Table1}.
\item Let $u$ be a unitary transvection. Then $C_{Y'}(u)$ acts
on
 $C_V(u)/[V,u]$ with orbits  of length $1$,  $6$ and  $9$.
\item If $u$ is a unitary transvection, $S_2 \le C_{Y}(u)$ has
order $3$ and $C_{C_V(u)/[V,u]}(S_2) \neq 0$, then $\dim
C_V(S_2)=2$.
\item \blue{ For $S \in \Syl_2(Y)$ and $S_1 = S \cap Y'$, every $S$-invariant subspace of $W$ of dimension at least $2$ contains $C_V(S_1)$. }
 \item $Y$ does not contain a fours group all of whose non-trivial
 elements are unitary transvections.
\item $C_V(u_4)$ is generated by non-singular vectors.
\item The $2$-rank of $Y$ is $4$.
\end{enumerate}
\end{proposition}
\begin{proof} (i) From \cite[page 27, page 47]{Atlas}, we see  that
$\Aut(\SU_4(2))$ and $\Sp_6(2)$ both possess exactly four conjugacy
classes of involutions.  \\

(ii) By Witt's lemma  $Y$ has exactly two orbits  on the non-zero
elements of  $V^\#$ and they correspond to the singular and the
non-singular vectors. Since $2^8-1$ does not divide $|X|$, these
orbits are also orbits under the action of $X$. Since the lengths of
the orbits are 135 and 120, using \cite[page 26, page 46]{Atlas} we get
the given  structure of the stabilizers.

(iii) As $Y$ contains a Sylow 3-subgroup of $X$, we find
representatives of all $X$-conjugacy classes of elements of order $3$ in
$Y$. By \cite[page 27]{Atlas} there are exactly three conjugacy
classes of elements of order 3 in $Y$, which we easily distinguish
by their action on $V$. We have elements, which are fixed point
free, which have centralizer of dimension $2$ and those which have
centralizer of  dimension 4. In particular, these elements are not fused in $X$.

Let $d \in Y$ have $2$-dimensional fixed space on $V$. Then as
$C_V(d)$ is perpendicular to $[V,d]$ we deduce that
$C_V(d)$ is non-singular (a $1$-dimensional non-singular
$\GF(4)$-space).

(iv) For the unitary transvection $u$ we have that $\dim [V,u] = 2$.
Suppose that $u$ is not a unitary transvection but $u \in Y'$. Then,
as $V$ supports the structure of a vector space over $\GF(4)$, we
have that $[V,u]$ is 2-dimensional and so $\dim [V,u] = 4$. \blue{ Suppose next that $u$
is an involution in $Y\setminus Y'$ and let $P$ be the stabilizer of a maximal isotropic space $W$ of $\GF(4)$-dimension 2 in $V$. Then $O_2(P)$ is elementary abelian of order 16 and $P/O_2(P) \cong \Sym(5) \cong \SU_2(4):2$.  Since $P$ contains a Sylow $2$-subgroup of $Y$, we may suppose that $u \in P$. Furthermore $W$ and $V/W$ are  natural $\SL_2(4)$-modules. As $u \not\in O_2(P) \le Y'$, we have that $\dim [W,u]=2 =\dim [ V/W,u]$. Hence we get that $\dim [V,u]\geq 4$ and so as $\dim V =8$, we have $\dim [V,u]=4$.}

(v) Let $u$ be a unitary transvection. Then $C_{Y'}(u)$ acts on
 $C_{V}(u)/[V,u]$ as the group $\GU_2(2) \cong
\Sym(3) \times 3$ and has three orbits one of length $1$, one of
length $6$ and one of length $9$.

(vi) From (v),  a Sylow 3-subgroup $S_1$ of $C_{Y'}(u)$ contains two
subgroups of order 3 whose centralizer in $C_{V}(u)/[V,u]$ is of
order 4 and two which are fixed point free. As the elements of order
three in $C_{Y'}(u)$ act the same way on $[V,u]$ as on $V/C_{V}(u)$,
the elements with fixed points on $C_{V}(u)/[V,u]$ have centralizer
in $V$ of dimension $2$, as by (iii) there are no elements of order
three which centralize a subspace of dimension $6$. Now by coprime action
we get that one subgroup of order three in $S_1$ centralizes in $V$
a subspace of dimension $4$ and acts fixed point freely on
$C_V(u)/[V,u]$, one acts fixed point freely $V$ and the other two
centralize a subspace of dimension $2$ in $V$.

(vii) \blue{Let $S\in \Syl_2(Y)$ and $S_1= S \cap Y'$. Then, as $V$ is the
natural $4$-dimensional unitary module for $Y'$, we have that
$U= C_V(S_1)$ has $\GF(2)$-dimension $2$. Assume that (vii) is false and let  $W$ be an $S$-invariant subspace of dimension at least $2$ with $U \not \le W$. Then $W>[W,S] \neq 0$ does not contain $U$ and so $[W,S]$ must have $\GF(2)$-dimension $1$ by the minimal choice of $W$.  Hence $[W,S]\le C_V(S_1)=U$ which means that $W +U/U\le C_{V/U}(S)< C_{V/U}(S_1)$ and this latter space  has $\GF(4)$-dimension $1$.  It follows that $W$ has $\GF(2)$-dimension $2$. Hence $S_0=C_{S}(W)$ has index $2$ in $S$,  $S_0 \cap S_1$ has order at least $2^5$ and this subgroup centralizes $W$ and $U$ and hence centralizes the preimage of $C_{V/U}(S_1)$ which has $\GF(4)$-dimension $2$. However, this is an isotropic line in the unitary representation and its centralizer is elementary abelian of order $2^4$, a contradiction. Hence (vii) is true.}

(viii) Suppose that $F = \langle x_1,x_2\rangle$ is a fours group
with all non-trivial elements unitary transvections. Then, as
$x_3=x_1x_2$, is also a unitary transvection, we get that
$C_V(x_1)= C_V(x_2)$. But then $C_V(x_1)$ is normalized by $\langle
C_Y(x_1), C_Y(x_2)\rangle = Y$, which is impossible.

(xi) Let $y$ be a non-singular vector. By (ii), we have that
$C_{Y}({y})\cong 3^{1+2}_+.\GL_2(3)$. This group contains an
involution $u$ in $Y\setminus Y'$. If $u$ is conjugate  to $u_3$ (in Table~\ref{Table1}), then
$C_{Y'}({u})\cong \Sym(6)$ acts transitively on $C_V(u)^\#$ and
so $C_V(u)^\#$ contains only non-singular vectors. Since $\dim
C_V(u)=4$, this is impossible. Therefore $v$ is conjugate to $u_4$
and $y \in C_V(u)=[V,u]$. Since $C_{C_{Y'}({y})}(u)$ has order
$6$, there are eight conjugates of $y$ in $C_V(u)$. Hence $C_V(u)$ is
generated by non-singular elements.

(x) \blue{ From (i) we see that the centralizers of involutions $x \in Y \setminus Y'$ have $2$-rank $4$. Thus we only need to see that $Y'$ has $2$-rank $4$. This is well-known and can be read from  \cite[Table 3.3.1]{GLS3}.}
\end{proof}

In the next lemma  the group denoted by $(\SU_4(2) \times 3){:}2$ is
the subgroup of index $2$ in $\Aut(\SU_4(2))\times \Sym(3)$ which is
not expressible as a direct product.

\begin{lemma}\label{notsimple} Assume that $G$ is a group, $t \in G$ is an involution, $H = C_G(t)$ and $Q=F^*(H)$ is
extraspecial of order $2^9$. If $H/Q \cong \Aut(\SU_4(2))$ or
$(\SU_4(2) \times 3){:}2$ and $Q/\langle t \rangle$ is the natural
$F^*(H/Q)$--module, then $G$ has a subgroup of index $2$.
\end{lemma}

\begin{proof}  We let $S \in \Syl_2(H)$ and note that, as $Z(S)=Z(Q) = \langle t \rangle$, we have $S \in \Syl_2(G)$.  Let $\ov H = H/\langle t\rangle$. We first show that

\bigskip
\begin{claim}\label{claim1}
$t^G\cap Q= \{t\}$.
\end{claim}

\bigskip

Assume that $u\sim_Gt$  with $u \in Q\setminus\langle t\rangle$.
Then $\ov u$ is singular in $\ov Q$ and so we may suppose that
$\langle \ov u \rangle = Z(\ov S)$.
 Now $C_Q(u)$ contains
an extraspecial group of order $2^7$.\blue{ As a Sylow 2-subgroup of
 $H/Q$ is not extraspecial}, we have that $t \in Q_u = O_2(C_G(u))$.  Note that
$\Phi(Q_u \cap Q) \le \langle u\rangle \cap \langle t \rangle =1$.
Hence $Q_u \cap Q$ is elementary abelian. As $Q$ is extraspecial of
order $2^9$, we deduce that $|Q\cap Q_u| \le 2^5$. Since the
$2$-rank of $H/Q$ is $4$ by \blue{Proposition~\ref{fact} (x)} and $|C_{Q_u}(t)|=2^8$, we infer that $|Q
\cap Q_u|$  is either $2^4$ or $2^5$. Furthermore, because $C_H(u)Q
\ge S$, we have that $Q \cap Q_u$ is a normal subgroup of $S$. We
know that $\ov Q$ is a $\GF(4)$-module for $F^\ast(H/Q)$.  Let
$\ov{U}$ be the one-dimensional $\GF(4)$-space in $\ov{Q}$
containing $\ov{u}$,  $U$ be its preimage in $H$ and set $R=
C_H(U)$. Since $U$, $Q_u\cap Q$ and $R$ are normalized by $S$,
Proposition~\ref{fact} (vii) implies $U \le Q_u \cap Q \le R$.
Assume that $|Q_u \cap Q|=2^5$.\blue{ Then, as $(Q_u \cap H) Q/Q$ is a
normal subgroup of $C_H(u)Q/Q$ and $C_H(u)Q/Q$ contains $S/Q$, we get $Z(S/Q) \leq (Q_u \cap H) Q/Q$. Hence  there exists $w\in Q_u\cap H$
such that $\langle wQ \rangle=Z(S/Q)$ is the unitary transvection
group centralizing $\ov R$.  Therefore we have
$$[Q_u\cap Q,w] \le [R,w] \cap [Q_u,w]\le \langle t \rangle \cap \langle u \rangle=1,$$} which is
impossible as $Q_u \cap Q$ is a maximal abelian subgroup of $Q_u$.
Thus $|Q_u \cap Q| = 2^4$. Since $|(Q_u \cap Q)/U|=2$, we now have
a contradiction to the fact that $C_{R/U}(C_H(u))=1$ by
Proposition~\ref{fact} (v). Thus \ref{claim1} holds.

\bigskip

 By Proposition~\ref{fact} (i), $H/Q$ has exactly two conjugacy classes of involutions not in $H^\prime/Q$. We choose representatives $\wt x, \wt y \in S/Q$ for these conjugacy classes and fix notation so that $C_{F^*(H/Q)}(\wt x) \cong \Sp_4(2)$ and $C_{F^*(H/Q)}(\wt y) \cong 2 \times \Sym(4)$.
We have that $|[\ov Q, \wt x]| = |[\ov Q, \wt y]| = 2^4$ by
Proposition~\ref{fact} (iv). Let $z \in H$  with $z^2 \in \langle t
\rangle$ be such that $zQ$ is either $\wt x$ or $\wt y$.  Let $T \in \Syl_2(C_H(z))$. Then $T' \cap
Z(T) \le T \cap H'$ and $Z(T) \cap H' \le Q $ as $Z(T) = \langle
z,C_Q(z)\rangle$. Thus, by \ref{claim1}, we have  $t^G \cap T' \cap
Z(T)= \{t\}$. In particular, $T\in \Syl_2(C_G(z))$.  It follows that
$z$ is not conjugate to $t$ in $G$ and that $t^G \cap Z(T)= \{t\}$.
We record these observations as follows:

\begin{claim}\label{claim2} Let $z \in S \setminus (S \cap H^\prime)$ be such that $z^2 \in \langle t\rangle$ and $T\in \Syl_2(C_H(z))$. Then $T\in \Syl_2(C_G(z))$, $t^G \cap Z(T) = \{t\}$ and $t^G \cap H \subset H^\prime$.
\end{claim}

\bigskip
Now let $z_1 \in S$ be such that $z_1Q= \wt x$.  Since
$C_{H/Q}(z_1Q)$  contains an element $fQ$ of order $5$ with $f$ of
order $5$  acting fixed point freely on $\ov Q$,  we see that $C_{
Q\langle z_1\rangle}(f)$ has order $4$. Let $z \in C_{Q\langle z_1 \rangle}(f)$ have
minimal order so that $zQ= z_1Q$. Then  $z^2 \in \langle t \rangle$.
Suppose that $g \in G$ and  $z^g\in S\cap H'$ is  extremal in $S$.
Then $C_S(z^g) \in \Syl_2(C_G(z^g))$. Now let $T \in
\Syl_2(C_H(z))$. Then $T \in \Syl_2(C_G(z))$ by \ref{claim2}. Hence
$T^g \in \Syl_2(C_G(z^g))$ and there is a $w \in C_G(z^g)$ such that
$T^{gw}= C_S(z^g)$. Now, by \ref{claim2}, $t^G \cap Z(T^{gw}) =
\{t^{gw}\}$ and of course $t^G \cap Z( C_S(z^g)) = \{t\}$ as $t \in
Z(H)$. Thus $gw \in H$, which is impossible as $z \in H \setminus
H'$, $z^g \in H'$ and $z^{gw} = z^g$. Hence there are no extremal
conjugates of $z$ in $S \cap H'$.  Since also $z^2 \in \langle t
\rangle$ and $t^G \cap H \subset H'$, Lemma~\ref{Thompson} implies
that $G$ has a subgroup of index $2$ as claimed.
\end{proof}

\section{The finer structure of $M$}

Suppose that $G$ is a group, $Z \le G$ has order $3$ and set $M = N_G(Z)$. \red{Assume that $C_M(Z)$ is similar
to a $3$-centralizer in a  group of type $\PSU_6(2)$ or $\F_4(2)$.} Let $S \in \syl_3(M)$ and $Q= F^*(M)=
O_3(M)$. By Hypothesis $C_M(Z)$ contains a normal subgroup $R^*$ such that $R^*/Q \cong \Q_8\times \Q_8$. We let
$R \in \Syl_2(R^*)$. Since the commutator map from $Q/Z\times Q/Z$ to $Z$ is an $C_M(Z)/Z$-invariant
non-degenerate symplectic form by \cite[III(13.7)]{Huppert} which may be negated by $M$, $M/Q$ embeds into
$\Out(Q) \cong \GSp_4(3)$. Our first lemma   locates $M/Q$ as a subgroup of $\GSp_4(3)$.

\red{\begin{lemma}\label{U6F4} We have that $M/Q$ normalizes $R^*/Q$ and is isomorphic to a subgroup of  the
subgroup $\mathbf M$  of $\GSp_4(3)$ which preserves a decomposition of the natural $4$-dimensional symplectic
space over $\GF(3)$ into a perpendicular sum of two non-degenerate $2$-spaces. Furthermore, $R/Q$ maps to
$O_2(\mathbf M)$.
\end{lemma}

\begin{proof}
Consider the action of $ Z(R)$ on $Q/Z$. Since $\Out(Q)$ is isomorphic to a subgroup of $\GSp_4(3)$, $Z(R)$ acts
as a fours group of $\Sp_4(3)$ on $Q/Z$.  Let $a \in Z(R)^\#$. Then $Q=C_Q(a)[Q,a]$  and $[C_Q(a),[Q,a]]=1$ by
the Three Subgroup Lemma.  We may suppose that  $C_Q(a) \neq Z$, and so we have $C_Q(a) \cong [Q,a]$ is
extraspecial of order $3^3$. Since $R$ centralizes $a$, $R$ preserves this decomposition and $R_1=C_R([Q,a])$ has
order $8$ and acts faithfully on $C_Q(a)$. Hence $R_1 \cong \Q_8$ and similarly $R_2 = C_R(C_Q(a)) \cong \Q_8$
with $R=R_1\times R_2$. In particular, we now have $C_M(Z) /Q$ is isomorphic to a subgroup of $\Sp_2(3) \wr 2$
and $R/Q$ corresponds to the largest normal $2$-subgroup of this group. It follows that $|O_2(C_M(Z)):RQ|\le 2$.
Thus $Z(R)Q/Q$ is a characteristic subgroup of $C_M(Z)$ and so $Z(R)Q/Q$ is normalized by $M/Q$. Finally, as
$RS/Q$ is the centralizer of   $Z(R)Q/Q $ in $C_M(Z)/Q$ we deduce that $RQ/Q$ is normalized by $M/Q$ and that
$M/Q$ preserves the decomposition of $Q/Z$ as described.
\end{proof}}

For the remainder of the paper we now assume  that $Z$ is weakly closed in $Q$ but not in $S$ with respect to
$G$. In particular, this means that $S > Q$.

\begin{lemma}\label{basic} The following hold.
\begin{enumerate}
\item $Z=Z(S)=Z(Q)$, $N_G(S)\le M$ and $S \in \syl_3(G)$;
\item $3\le  |S/Q|\le 3^2$; and
\item $Q$ has exponent $3$.\end{enumerate}
\end{lemma}

\begin{proof}

(i) Since $C_M(Q) \le Q$, we have that $Z= Z(Q)=Z(S)$.  Therefore
$N_G(S) \le N_G(Z)= M$ and, in particular, $S
\in\Syl_3(N_G(S))\subseteq \Syl_3(G)$.

 (ii) This follows directly from Lemma~\ref{U6F4}.

(iii) Since $[Q,a]$ admits $R$ and $C_Q(a)$ admits $R$, these groups have exponent $3$ and they commute. Thus
(iii) holds.
\end{proof}

 We  have  $\ov M = M/Q$ is isomorphic to a subgroup of  the subgroup of $\GSp_4(3)$ which preserves a
decomposition of the natural $4$-dimensional symplectic space into a perpendicular sum of two non-degenerate
$2$-spaces by Lemma~\ref{U6F4}. We now describe this subgroup of $\GSp_4(3)$. We denote it by $\ov{{\mathbf M}}$
as in Lemma~\ref{U6F4}. The boldface type is supposed to indicate that this is a subgroup of $\GSp_4(3)$ which
contains (the image of) $\ov M$ but may be greater than it. Similarly $\ov {\mathbf S}$ is a Sylow $3$-subgroup
of $\ov{{\mathbf M}}$ which contains $\ov S$.

We have $\ov {\mathbf M}$ contains a subgroup of index $2$ which is
contained in $\Sp_4(3)$ and is isomorphic to the wreath product of
$\Sp_2(3) \cong \SL_2(3)$ by a group of order $2$. For $i=1, 2$, we
let $\ov{\mathbf {M_i}}\cong \SL_2(3)$,  $\ov{\mathbf {R_i}} =
O_2(\ov {\mathbf M_i})\cong \Q_8$ and $\ov{\mathbf {S_i}} =
\ov{\mathbf S} \cap \ov{\mathbf{M_i}}$. We let $\ov{{\mathbf t_1}}$
be an involution in  $\ov {{\mathbf M}}$ which negates the
symplectic form and normalizes $\ov {{\mathbf S_1}}$ and $\ov
{{\mathbf S_2}}$. Note that, for $i=1,2$, $\ov{\mathbf M_{i}}\langle
\ov{{\mathbf t_1}}\rangle \cong \GSp_2(3)\cong \GL_2(3)$.  Next
select an involution $\ov {{\mathbf t_2}}$ which commutes with
$\ov{{\mathbf t_1}}$, preserves the symplectic form, normalizes $\ov
{{\mathbf S}}$ and conjugates $\ov{\mathbf {M_1}}$ to $\ov{\mathbf
{M_2}}$. With this notation we  have
$$\ov{{\mathbf M}}=\ov{\mathbf {M_1}}\ov{\mathbf {M_2}}\langle  \ov {{\mathbf t_1}},\ov {{\mathbf
t_2}}\rangle.$$

Now $\ov{M}$ is a subgroup of $\ov{{\mathbf M}}$ which has index at most $6$. In particular, $\ov S$ has index at
most $3$ in   $\ov {{\mathbf S}}$ by Lemma~\ref{basic} (ii). Since $\ov {{\mathbf R_1}}  \ov {{\mathbf R_2}} =\ov
R$, $M$ contains subgroups $R_1$ and $R_2$ isomorphic to $\Q_8$ such that $[R_1,R_2]=1$ and $\ov{R_i} =
\ov{{\mathbf R_i}}$ for $i=1$ and $2$.  Moreover $R= R_1R_2$. Let $T \in \syl_2(M)$ with $T \ge R$.   Now we do
not yet know the index of $R$ in $T$. Thus $T$ may or may not contain elements which map to $\ov{{\mathbf t_2}}$,
$\ov{{\mathbf t_2}}$ or $\ov{{\mathbf {t_1t_2}}}$. However if such elements are contained in $T$ we denote this
involution by $t_1$, $t_2$ or $t_1t_2$ as appropriate.   This discussion proves the following lemma.

\red{
\begin{lemma}\label{Syl2M}
There are exactly five possibilities for a Sylow $2$-subgroup $T$
of $M$. Moreover, one of the following hold.
\begin{enumerate}
\item $T=R$, $N_M(S) = SZ(R)$ and $N_M(S)/S \cong 2^2$;
\item $T= R\langle t_1\rangle$, $N_M(S) = SZ(R)\langle t_1\rangle$ and $N_M(S)/S \cong 2^3$;
\item $T= R\langle t_2\rangle$, $N_M(S) = SZ(R)\langle t_2\rangle$ and $N_M(S)/S \cong \Dih(8)$;
\item $T= R\langle t_1t_2\rangle$, $N_M(S) = SZ(R)\langle t_1 t_2\rangle$ and $N_M(S)/S \cong \Dih(8)$;
 and
\item $T=R\langle t_1,t_2\rangle$, $N_M(S) = SZ(R)\langle t_1,t_2\rangle$ and $N_M(S)/S \cong 2 \times \Dih(8)$.\qed\end{enumerate}
\end{lemma}

}

For $i=1,2$, let $r_i \in Z(R_i)^\#$ and set $Q_i =
[Q,r_i]=[Q,R_i]$.
Note that, as $\ov {r_1r_2} \in Z(\ov M)$ and $Q/Z$ is irreducible
as an $\ov M$-module, $r_1r_2$ inverts $Q/Z$. Let $A $ be the
preimage of $C_{Q/Z}(S)$. So $A$ is the second centre of $S$.

\begin{lemma}\label{Qaction} The following hold.
\begin{enumerate}
\item $Q_1=[Q,R_1]= C_Q(R_2)$, $Q_2= [Q,R_2]=C_Q(R_1)$ and both are normal in $S$;
\item $Q_1 \cong Q_2 \cong 3^{1+2}_+$, $[Q_1,Q_2]=1$ and $Q=Q_1Q_2$;
\item $A=[Q,S] = [Q_1,S][Q_2,S]$ is elementary abelian of order $3^3$; and
\end{enumerate}
\end{lemma}

\begin{proof}
(i) This follows directly from the action of $M$ on $Q$ as $\ov
{R_1}$ and $\ov {R_2}$ are normalized by $\ov S$.

(ii) We have that $C_Q(r_1)$ and $[Q,r_1]$ commute by the Three
Subgroup Lemma. Since, for $i=1,2$,  $[Q,r_i]=[Q,R_i]$ has order
$3^3$ it follows  that $Q_1\cong 3^{1+2}_+$. As $r_1r_2$ inverts
$Q/Z$, $r_2 $ inverts $C_{Q/Z}(r_1)$ and so  $C_Q(r_1)= Q_2$. In
particular, $Q_1$ and $Q_2$ commute and $Q=Q_1Q_2$.

(iii) From the description of $M/Q$, we have $A=[Q_1,S][Q_2,S]$.
\blue{Since, for $i = 1,2$, $[Q_i/Z,S] \not= 1$ and $[Q_i,S]$ is normal in $Q_i$,  we have $Z \leq [Q_i,S]$ and so  $[Q_i,S]$ has order $9$.} Furthermore $[Q_i,S]$ is elementary
abelian. Hence $A$ is elementary abelian of order $3^3$ by (ii).
\end{proof}

Because, for $i=1,2$, $r_i$ inverts $Q_i/Z$, if  $M$ happens to
contain the involution $t_1$, we may and do adjust $t_1$ by
multiplying by elements from $Z(R)$ so that $t_1$ inverts $A/Z$.
Therefore
\red{\begin{lemma}\label{ift1} If $t_1 \in M$, then  $t_1$ inverts $A$ and centralizes $Q/A$.\qed\end{lemma}}

We now define a subgroup which will play a prominent role in all the future investigations. Set  $$J= C_S(A).$$ It will turn out that $J$ is the Thompson subgroup of $S$.

\begin{lemma}\label{NMJ} The following hold:
\begin{enumerate}
\item $|S:J|=3^2$, $J\cap Q= A$ and $S=JQ$;
\red{\item $N_M(S) = N_M(J)$; and
\item if $t_1 \in M$, then $t_1$ inverts $J$ and $J$ is abelian.}
\end{enumerate}
In particular, the structure of  $N_M(J)/S$ is as described in the five parts of Lemma~\ref{Syl2M}.
\end{lemma}

\begin{proof} By Lemma~\ref{Qaction}(iii), $A$  is elementary abelian of order $3^3$. Furthermore, by the definition of $J$, $J $ is a normal subgroup of $N_M(S)$.  Since $[S,A]= Z$,
the $3$-structure of $\GL_3(3)$ shows that $|S/J|\le 3^2$. As $J
\cap Q= C_Q(A)=A$, we infer that $|S:J|= 3^2$ and $S=JQ$.  Thus (i)
holds.

\red{As $A = Z_2(S)$,  we have $J=C_S(A)$ is normalized by $N_M(S)$. Also because $N_M(J)$ normalizes $J \cap Q=A$, we know  $N_M(J)\le  N_M(A)$.
 Now
  $N_M(A) \ge N_M(S)$ and because $N_R(A) = Z(R)$ and $M= RN_M(S)$,  we deduce that $N_M(A)= N_M(S)$.  Therefore $N_M(J)= N_M(S)$ as claimed in (ii).}

Suppose that  $t_1 \in M$.
Then  $t_1$ inverts $S/Q$, centralizes $Q/A$ and inverts $A$ by
Lemma~\ref{ift1}.
 Thus $t_ 1$ inverts $J$ and so $J$ is
abelian. This concludes the proof of (ii) and completes the
verification of the lemma.
\end{proof}

Note that $|J|=3^4$ if $|S/Q|=3$ and $|J|=3^5$ if $|S/J|=3^2$.

\begin{lemma}\label{CSJ} We have $C_G(J)= J$.
\end{lemma}

\begin{proof}\red{As $Z \le J$, we have $C_G(J) = C_M(J)$.  Now $C_G(J)$ centralizes $A=J \cap Q$ and it follows from Lemma~\ref{w=1} that $C_M(J) = C_S(J) = J$.}
\end{proof}

\begin{lemma}\label{ActionQ} Every element of $Q$ is conjugate in $M$ to an element of $A$.
\end{lemma}

\begin{proof} It suffices to prove that every element of $Q/Z$ is conjugate to an element
of $A/Z$. Let $w\in Q/Z$. Then $w=x_1x_2$ where $x_i \in Q_i/Z$ by
Lemma~\ref{Qaction} (ii). Since, by Lemma \ref{Qaction} (iii), for
$i=1, 2$, $(A\cap Q_i)/Z$ has order $3$ and $R_i$ acts transitively
on $Q_i/Z$, there exists $u_i\in R_i$ such that $w^{u_1u_2} =
x_1^{u_1}x_2^{u_2} \in A/Z$. This proves the claim.
\end{proof}

\section{The structure of the normalizer in $G$ of $J$}

For the remainder of the paper  assume the hypothesis of the
Theorem~\ref{MainTheorem}. Thus we have $M,Q, S$ and $Z$ as in
Section~3 and additionally we have that $Z$ is weakly closed in $Q$
and not in $M$. In this section we determine the possible structures
of $N_G(J)$.

\begin{lemma}\label{J ab} If $Z$ is not weakly closed in $J$, then $J$ is elementary abelian and coincides with the Thompson subgroup of
$S$. In particular, $N_G(J)$ controls fusion in $J$.
\end{lemma}

\begin{proof} Choose $X \in Z^G$ with $X\neq Z$ and $X \le J$. Set $K= AX$. As $Z$ is weakly closed in $Q$ and $J = C_S(A)$, we have that $K$ is elementary abelian of order $3^4$. In particular, if  $|J|= 3^4$, then $K = J$ is elementary abelian. \\
\\
Suppose that $|J|=3^5$.  Then $|J : K| = 3$ and $|S/Q|=3^2$. We claim
that $J$ is abelian.
   Set $Q_X = O_3(N_G(X))$.
As $K$ has index $3$ in $J$, $K$ is normal in $J$ and, as $[Q,X] \le
A$, $K$ is normalized by $Q$. Therefore $K$ is normal in $ S = JQ$
by Lemma~\ref{NMJ} (i). If $C_S(X)= K$, then $|X^S|= 3^3$ and, in
particular, every element of $K$ which is not conjugate to an
element of $Z$ is contained in $A$. Now $K \cap Q_X$ has order
either $3^2$ or $3^3$ and, so, as  $X$ is weakly closed in $Q_X$,
$K\cap Q_X$ is generated by elements which are not conjugate to
elements of $Z$. It follows that $X\le K\cap Q_X \le A$ and this
contradicts $X\not \le Q$. Therefore $C_S(X) \neq K$. If
$C_S(X)\not\le J$, then
 $Z= [A,C_S(X)] \le C_S(X)'\le Q_X$ and this contradicts the fact that $X$ is
weakly closed in $Q_X$. So  $C_S(X)\le J$. But  then we have $K \le
Z(J)$ and so $J$ is abelian as claimed.

Suppose that $B \le S$ is abelian and $|B| \ge |J|$. Then, as
$|B\cap Q| \le 3^3$, we have $BQ=S$ and then $(B\cap Q)/Z \le
C_{Q/Z}(S) = A/Z$. Thus $B \le C_S(A) = J$. Hence $J$ is the
Thompson subgroup of $S$. \blue{As $J$ is abelian and  weakly closed in $S$, it} follows from \blue{\cite[(37.6)]{Aschbacher}} that $N_M(J)$ controls fusion
in $J$.  In particular, $X$ and $Z$ are conjugate in $N_M(J)$. Since
$\Phi(J) \le A$, $X \not \le \Phi(J)$ and hence $Z \not \le
\Phi(J)$.  Therefore $Z(S) \cap \Phi(J)=1$. As $\Phi(J)$ is normal
in $S$, we get $\Phi(J) = 1$  and $J$ is elementary abelian. This
completes the proof of the lemma.
\end{proof}

\begin{lemma}\label{Z 10} Assume that $Z$ is not weakly closed in $J$ and set $J_0 =\langle Z^{N_G(J)} \rangle$.
Then
  \begin{enumerate}\item $|Z^{N_G(J)}| = 10$ and, if $X \in Z^{N_G(J)}$ with $X \neq Z$, $|X^Q|=3^2$;
   \item $N_G(J)$ acts $2$-transitively on $Z^{N_G(J)}$; and
   \item $|J_0 Q/Q| =3$ and $J_0 Q/Q$ is normalized by $N_M(S)/Q$.\end{enumerate}
\end{lemma}

\begin{proof} Let $\mathcal Y =Z^{N_G(J)}$ and $X \in \mathcal Y$ with $X \neq Z$. Of course $X \not \le Q$ as $Z$ is weakly closed in $Q$. If $C_Q(X) \not \le J$, then, as $X$ centralizes $A$,
$C_Q(X)$ has order at least $3^4$ and consequently is non-abelian and we have
$Z=C_Q(X)'\le O_3(C_G(X))$. However $X$ is weakly closed in
$O_3(C_G(X))$ with respect to $G$ and so this is impossible. Thus
$C_Q(X)= A$ has order $3^3$ and, in particular (as $J$ is abelian), $X^S= X^{JQ} = X^Q$
has order $3^2$ and so $|\mathcal Y|\equiv 1\pmod 9$. \blue{Observe that
 $$|\mathcal Y| = |N_G(J)|/|N_M(J)| =  |N_G(J)/J| /|N_M(J)/J|.$$} As $|J|= 3^4$ or $3^5$ and $J$ is
self-centralizing and elementary abelian by Lemmas~\ref{CSJ} and
\ref{J ab}, $|N_G(J)/J|$ divides $|\GL_5(3)|$. If $|J|= 3^4$, then,
as no subgroup of order $3$ in $A$ which is not $Z$ is conjugate
to $Z$, $J$ contains at most $28$ conjugates of $Z$. This means that
$|\mathcal Y|= 10, 19$ or $28$. On the other hand,
$|\GL_4(3)|_{3^\prime}= 2^9\cdot 5\cdot 13$ and so in this case
$|\mathcal Y| = 10$. So assume from now on that $|J|= 3^5$.  Then
$J$ contains $121$ subgroups of order $3$ and $12$ of these are
contained in $A$ and  are not conjugate to $Z$ as $Z$ is weakly
closed in $Q$. Since $|\GL_5(3)|_{3^\prime}= 2^{10}\cdot 5\cdot
11^2\cdot 13$ and $|\mathcal Y| \equiv 1 \pmod 9$, the only
candidates for $|\mathcal Y|$ are $10$, $55$ and $64$. We recall
from Lemma~\ref{NMJ} that $|N_M(J)/J| = 2^{i}\cdot 3^2$ where $i \in \{2,3, 4\}$ and, if $t_1 \in M$, then $t_1J \in
Z(N_G(J)/J)$  and $t_1$ inverts every element of $J$ by Lemma~\ref{NMJ}(iii).  In particular, $t_1$ normalizes every
member of $\mathcal Y$.

Suppose that $|\mathcal Y|= 55$. Then, by Lemma~\ref{NMJ} (ii),  $$|N_G(J)/J| = |N_G(J):N_{M}(J)||N_M(J):J|=2^i\cdot
3^2\cdot 5\cdot 11$$ where $i \in \{2,3,4\}$. Let $E \in
\syl_{11}(N_G(J)/J)$. Then, as the normalizer of a cyclic subgroup
of order $11$ in $\GL_5(3)$ has order $2\cdot 5\cdot 11^2$, the
normalizer in $N_G(J)/J$ of $E$ has order dividing $110$. In
particular, $E$ is not normal in $N_G(J)/J$. If $|N_M(J)|_2=2^4$,
then $t_1J$ normalizes $E$.  So in any case the number of conjugates
of $E$ in $N_G(J)/J$ divides $2^3\cdot 3^2\cdot 5$ and is divisible
by $2^2\cdot 3^2$ and this is  impossible as it must also be
congruent to $1$ mod $11$.

Suppose that $|\mathcal Y|= 64$. Then $|N_G(J)/J| = 2^{j}\cdot 3^2$
where $j \in \{8,9,10\}$.  In particular, $N_G(J)$ is soluble.  Since
$|\mathcal Y|= 64$, we have that $J = \langle \mathcal Y\rangle$. If
$1\neq K \le J$ is normal in $N_G(J)$, then $K$ is normal in $S$ and
consequently $Z \le K$. But then $\mathcal Y\subseteq K$ and so
$K= J$. Thus $N_G(J)$ acts irreducibly on $J$. Since $|J| = 3^5$ and
$N_G(J)/J$ is not abelian, Schur's Lemma implies that
$|Z(N_G(J)/J)|$ divides $2$ and, additionally, $O_3(N_G(J)/J) =1$.
Let $L = O_{3,2}(N_G(J))$. By Clifford's Theorem \cite[Theorem
4.3.1]{Gorenstein}, $J$ is completely reducible as an $L$-module
and $N_G(J)$ acts transitively on the homogeneous summands of $J$
restricted to $L$. Since $J$ has dimension $5$ as a
$\GF(3)N_G(J)$-module, and $5$ does not divide $|N_G(J)|$, we have
that $J$ is homogeneous as an $L$-module. It follows that $J$ is
either  a direct sum of five $1$-dimensional $L$-modules or is
irreducible as an $L$-module. It the first case, we get that
$[L,N_G(J)]\le J$, $O_3(N_L(J)/J) \not= 1$ and this contradicts $O_3(N_G(J))=J$. Thus $J$ is
an irreducible $L$-module. However, the degrees of irreducible
$L/Q$-modules over the algebraic closure of $\GF(3)$ are all powers
of $2$ \cite[15.13]{Isaacs} and this again implies that $L/C_L(J)$ is
cyclic, and  $O_3(N_G(J)) > J$, again a contradiction. Thus $|\mathcal Y| \neq 64$.

Since $|\mathcal Y|\not = 55$ or
$64$, we must have $|\mathcal Y|=10$ as claimed in the first part of
(i). Because  we have also shown that $C_Q(X)=A$ the remaining parts of
(i) also hold.

Part (ii) follows directly from (i).

Now with $J_0 =\langle Z^{N_G(J)} \rangle$, we have  $\langle X^Q
\rangle Q= XQ$ is normalized by $N_M(S)$ and $|XQ/Q| = 3$. This is
(iii).
\end{proof}

\begin{lemma}\label{diagonal} Suppose that $X \in Z^G\setminus \{Z\}$  and $X \le S$. Then, for $i=1,2$, $[X,R_i]\not
\le Q$.
\end{lemma}

\begin{proof} We  suppose that $[X,R_1] \le Q$ and seek a contradiction. Let $Q_X =
O_3(N_G(X))$ and $W$ be the full preimage of $C_{Q/Z}(X)$. Since $R_1$ acts irreducibly on $Q_1/Z$ and $[Q_1,QX]$ is
$R_1$-invariant, we have $Q_1 \le W$.  Then
$|W|=3^4$. Hence $W= Q_1A\cong 3 \times
3^{1+2}_+$ and $Z(W)= A\cap Q_2$.

 If $C_W(X)$ is
non-abelian, then, as $C_W(X)Q_X/Q_X$ is abelian, $Z=C_W(X)'\le
Q_X$. Since $X$ is weakly closed in $Q_X$ by assumption and $Z \neq
X$, we have a contradiction.

 Thus $C_W(X)$ is abelian. Since $W$ is
non-abelian and $XZ$ is normalized by $W$, we get that
$|C_W(X)|=3^3$.  Because $C_W(X)$ is abelian and $W$ is not, it
follows that $A \cap Q_2 \le C_W(X)$. Furthermore, we have $|C_W(X)
\cap Q_1|= 3^2$ and thus, as $R_1$ acts transitively on the
subgroups of order $9$ in $Q_1$, we may adjust $X$ by conjugating by
an element of $R_1$ and arrange for $W \cap Q_1 = A \cap Q_1$. But
then $W = A$ and $X \le J$. Put $J_0 = \langle X^{N_G(J)}\rangle$.
Then by Lemma~\ref{Z 10} (iii)  $J_0Q = XQ$ is normalized by
$N_M(S)$. Since $N_M(S)$ does not normalize $R_1$, we have
$[X,R_1R_2] \le Q$, and this contradicts the structure of $M$.
Therefore $[X,R_i] \not \le Q$ for both $i=1$ and $2$.
\end{proof}

\begin{lemma}\label{Z in J}
Assume that $X\in Z^G$ with $X \le S$. Then $X \le J$. In
particular, $Z$ is not weakly closed in $J$.
\end{lemma}

\begin{proof} Suppose that $X \le S$ and $X \not \le J$. Then $[A,X]= Z$ and $|C_A(X)|= 3^2$. By
Lemma~\ref{diagonal}, $XQ$ acts non-trivially on both $R_1Q/Q$ and
$R_2Q/Q$ and so $C_A(X) = C_Q(X)$. On the other hand $AX$ is
normalized by $Q$ and so $AX$ contains at least, and hence exactly,
$28$ conjugates of $Z$. In particular, $C_A(X)X$ contains $10$
conjugates of $Z$ and three subgroups of order $3$ which are not
conjugate to $Z$. Set $Q_X= O_3(N_G(X))$. Then the only conjugate of
$Z$ contained in $C_A(X)X \cap Q_X$ is $X$. Since the subgroups of
order $3$ in  $C_A(X)$ which are not conjugate to $Z$ generate
$C_A(X)$, we get $C_A(X)X \cap Q_X = X$. So $|C_A(X)Q_X/Q_X|= 3^2$.
By Lemma~\ref{diagonal} two of the non-trivial cyclic subgroups of
$C_A(X)Q_X/Q_X$ are not images of elements  from $Z^G$. Since
$C_A(X)X$ contains only three subgroups of order $3$  which are not
conjugate to $Z$, we have a contradiction. Therefore, if $X \in Z^G$
and $X \le S$, $X \le J$ as claimed.
\end{proof}

 Set $$J_0 = \langle Z^{N_G(J)}\rangle.$$ By Lemmas~\ref{Z 10}, \ref{diagonal} and \ref{Z in J}, we have $|J_0 Q/Q| = 3$, $J_0 \cap Q = A$ and ${J_0}Q/Q$ does not centralize either ${R_1}Q/Q$ or ${R_2}Q/Q$. In particular, $|J_0| = 3^4$. We  record these facts in the first part of the next lemma.

\begin{lemma}\label{NGJ0} The following hold.
\begin{enumerate}
\item $|J_0| =3 ^4$, $|J_0Q/Q|=3$, $J_0 \cap Q=A$ and $J_0Q/Q$ acts non-trivially on both $R_1Q/Q$ and $R_2Q/Q$;
\item $N_G(J)= N_G(J_0)$; and
\item $C_G(J_0)= C_G(J)=J$.

\end{enumerate}
\end{lemma}

\begin{proof} From the construction of $J_0$ we have $N_G(J_0)\ge N_G(J)$. Since $N_G(J)$ is transitive on the subgroups of $J$ which are $G$-conjugate to $Z$, we get that $N_G(J_0) = N_G(J)N_M(J_0)$. Hence, as $N_M(J_0Q) = N_M(S) \le N_G(J)$,  (ii) holds. Obviously $C_G(J_0) \le C_M(J_0) \le C_M(A)= J$ so (iii) also holds.
\end{proof}

Define $$F = O^2(N_G(J))\langle r_2\rangle.$$ Note that $F$ is a group as $r_2$ normalizes $A/Z$ and hence $J$. Then

\begin{theorem}\label{S6}The following hold:
\begin{enumerate}
\item The action of $N_G(J)$ on $J_0$ preserves a non-degenerate quadratic form $\mathrm q$ of $(-)$-type;
\item $Z^{N_G(J)}$ is the set of singular one-dimensional subspaces with respect to $\mathrm q$;
\item $N_G(J)/J\cong 2\times \Sym(6)$ or $\Sym(6)$; and
\item $F/J \cong \Sym(6)$ and $|[J,r_2]|=3$. Furthermore $[r_2,J]\le J_0$ and $[J,F] \le J_0$.
\end{enumerate}
\end{theorem}

\begin{proof} Let $X \le J$ be conjugate to $Z$ but not equal $Z$. For $i=1,2$, using Lemma~\ref{NGJ0} (i), we have
that \blue{ $[J_0,Q_i/Z] \not= 1$ and so, as $[J_0,Q_i]$ is normal in $Q_i$, we get}  $|[J_0,Q_i]|= 3^2$ and $[J_0,Q_i,Q_i]= Z$. Furthermore,
$[J_0,Q_i]$ is centralized by $Q_{3-i}$. Hence we have
$[J_0,Q_{i}]=C_{J_0}(Q_{3-i})$. By Lemma~\ref{GO4}, there exists a
non-degenerate quadratic form $\mathrm q$ on $J_0$ which is
preserved by $Q$ and such that the elements of $X$ are singular
vectors. It follows that with respect to $\mathrm q$, the elements
of $\bigcup X^Q$ are singular. Furthermore, as $Z= C_{J_0}(Q)$, $Z$
also consists of singular vectors. Now with respect to the
bilinear form $\mathrm f$ associated with $\mathrm q$, none of the
non-trivial elements of $\bigcup X^Q$ are perpendicular to the
non-trivial elements of $Z$. It follows that $XZ$ contains exactly
two singular subspace, namely $X$ and $Z$. Since $N_G(J)$ acts $2$-transitively on $Z^{N_G(J)}$ by Lemma~\ref{Z 10} (ii), we infer that
if $X, Y\in Z^{N_G(J)}$ with $X \neq Y$, then $XY$ contains exactly
two members of $Z^{N_G(J)}$. Now suppose that $a \in Q\setminus J$
is such that $aJ$ acts quadratically on $J_0$.\blue{For $X \in
Z^{N_G(J)}\setminus Z$ we know $[X,a] \not= 1$ and hence $|[X,a]|= 3$ as $a$ acts quadratically on $J_0$.} It
follows that $X[X,a] = XX^{a}$ contains three members of $Z^{N_G(J)}$ namely
$X$, $X^a$ and $X^{a^2}$. This contradiction shows that no
non-trivial element of $S/J$ acts quadratically on $J_0$. If
$\mathrm q$ was of $(+)$-type, this would not be the case. Hence
$\mathrm q$ is of $(-)$-type. We now have that $Z^{N_G(J)}$ is the set
of singular one spaces in $J_0$ with respect to $\mathrm q$. Since
$N_G(J)$ preserves this set, we have that $N_G(J)/J$ is isomorphic
to a subgroup of \red{$\mathrm {CO}_4^-(3)$, the group preserving $\mathrm q$ up to negation, from Lemma~\ref{quadratic
form}. Because $N_{N_G(J)}(Z)$ has index $10$ in $N_G(J)$, we deduce
that $|N_G(J)|= 2^i\cdot 5 \cdot 3^2$  where $i \in \{2,3,4\}$  by  Lemma~\ref{Syl2M}.
 In particular, $O^2(N_G(J)/J) \cong \Omega_4^-(3)\cong \Alt(6)$.  Now we see that $N_{O^2(N_G(J))}(S)/S$ is a cyclic group of order $4$. Consequently of the five possibilities for the structure of $N_M(S)/S$ given in Lemma~\ref{Syl2M} only possibilities (iii), (iv) and (v) survive and we have $N_M(S)/N_{O^2(N_G(J))}(S)$ is elementary abelian.

 Now let $C = \mathrm {CO}_4^-(3)$, if $D \in \Syl_3(X)$ and $Y \in \Syl_2(N_C(D))$, then $Y \cong 2 \times \SDih(16)$
and so using the structure of $N_M(S)/S$ given in Lemma~\ref{Syl2M} (iii), (iv) and (v) we infer
that $N_G(J)/J \cong \Sym(6)$ or $\GO_4^-(3) \cong 2\times \Sym(6)$.
We have now established (i), (ii) and (iii).}

We know that $S/Q = JQ/Q$ is centralized by $r_2$ and that $[Q,r_2]=
Q_2$. It follows that $[J,r_2] \le Q_2\cap J = A \cap Q_2$ and, as
$[A\cap Q_2,r_2]$ has order $3$, we now have that $[J,r_2]= [A,r_2]$
is a non-central cyclic subgroup of $Q$. In particular, $[J,r_2] \le
A\le J_0$. Since $|[J_0,r_2]| = 3$ we get that $r_2$ has determinant
$-1$ on $J_0$. Hence we have $r_2\not \in O^2(N_G(J))$ and so we
conclude that $F/J\cong \Sym(6)$ and that all the parts of (iv)
hold.
\end{proof}

As a corollary to Theorem~\ref{S6} we record the following observation.

\begin{corollary}\label{cor} There are exactly three possibilities for a Sylow $2$-subgroup $T$
of $M$. They are as follows:
\begin{enumerate}
\item $T= R\langle t_2\rangle$, $N_M(S) = SZ(R)\langle t_2\rangle$ and $N_M(S)/S \cong \Dih(8)$;
\item $T= R\langle t_1t_2\rangle$, $N_M(S) = SZ(R)\langle t_1 t_2\rangle$ and $N_M(S)/S \cong \Dih(8)$;
 and
\item $T=R\langle t_1,t_2\rangle$, $N_M(S) = SZ(R)\langle t_1,t_2\rangle$ and $N_M(S)/S \cong 2 \times \Dih(8)$.\qed\end{enumerate}
In particular, $Q/Z$ is a chief factor in $M$.
\end{corollary}

\begin{proof} The first three statement are readily deduced from the structure of $N_G(J)/J$ and so we only need to
elucidate the fact that $Q/Z$ is a chief factor. For this we simply note that $t_2 \in M$ in all cases.
\end{proof}

\begin{theorem}\label{NMSagain}
If  $N_M(S)/S\cong \Dih(8)$, then $G \cong \PSU_6(2)$ or
$\PSU_6(2){:}3$.
\end{theorem}

\begin{proof} Since $N_M(S)/S \cong \Dih(8)$, we have that $N_G(J)/J \cong \Sym(6)$ from the proof of Theorem~\ref{S6} (ii) and Corollary \ref{cor}. Since $[J,r_2]\le J_0$, we infer that $J/J_0$ is centralized by $N_M(J)$.
If $J > J_0$, then,  by Lemma~\ref{ptransfer}, $G$ has a normal
subgroup $G^*$ of index $3$. If $J = J_0$, then set $G = G^\ast$.
Now $M \cap G^*$ satisfies the hypothesis of Theorem~\ref{Parker}.
Hence $G^* \cong \PSU_6(2)$ and this proves the theorem.
\end{proof}

In light of Theorem~\ref{NMSagain} and Corollary \ref{cor},  from here
on we may assume that $N_M(S) = SZ(R)\langle t_1,t_2\rangle$. In
particular from Theorem \ref{S6}, we have
\begin{eqnarray*}N_M(S)/S&\cong& 2\times \Dih(8);\\N_G(J)/J &\cong& 2\times \Sym(6); \text { and} \\C_{F/J}(r_2J) &\cong& 2 \times \Sym(4).\end{eqnarray*}
Furthermore, as $t_1$ inverts $J$, we have $t_1 J \in Z(N_G(J)/J)$.

\begin{lemma}\label{CQ1} We have $$C_S(Q_1)= C_S(R_1)= C_S(Q_1R_1),$$ $$C_J(Q_1)= C_J(R_1)=C_J(Q_1R_1)$$ and $|J:C_J(Q_1)|=3^2$.
\end{lemma}

\begin{proof} We have that $[Q_1,C_S(R_1)]$ is $R_1$-invariant and is a proper subgroup of $Q_1$. Therefore $[Q_1,C_S(R_1)] \le Z$. Hence $[Q_1,C_S(R_1),R_1]=1$ and $[C_S(R_1),R_1,Q_1]=1$ and thus the Three Subgroups Lemma implies that $[Q_1,R_1,C_S(R_1)]=1$. Since $Q_1=[Q_1,R_1]$, we have $C_S(R_1) \le C_S(Q_1)$. Now, as $Q_1$ is normal in $S$ and extraspecial of order $3^3$,  $|S:C_S(Q_1)Q_1| =3$, and so $|C_S(Q_1)| = 3^4$ if $|S| = 3^7$ and $|C_S(Q_1)| = 3^3$ if $|S|= 3^6$. Since $R_1$ centralizes $Q_2$, we have $C_S(R_1)= C_S(Q_1)= Q_2$ if $|S| = 3^6$. If $|S|= 3^7$,
then, as $R_1Q$ is normalized by $R_1S$, we have $|S/C_S(R_1)Q| =3$
and hence the equality $C_S(Q_1)=C_S(R_1)$ holds in this case as well. Of
course we now have  $C_J(Q_1)= C_J(R_1)= C_J(Q_1R_1)$.

Since $J$ normalizes $R_1Q$ and does not centralize $R_1Q/Q$ by
Lemma~\ref{diagonal}, $Q_1$  is normalized by $J$. Since $J$ is
abelian and $J \cap Q_1 = A\cap Q_1$, we now have that
$|J:C_J(Q_1)|=3^2$.
\end{proof}

\blue{Notice that $r_1J$ and $r_2 J$ are conjugate in $N_G(J)/J$ (by
$t_2J$ for example) and  $$\langle r_1, r_2 ,Q_1\rangle J/J \cong
2\times \Sym(3).$$ In particular, we have $r_1 \in F$.

Let $U\le F$ be chosen so that $\langle r_1, r_2 ,Q_1\rangle J \le U
$ and $U/J \cong \Sym(5)$.

\begin{lemma}\label{6orb} If $J \not=J_0$, then $|C_J(U)|=3$ and $|C_J(U)^F|= |C_J(U)^{N_G(J)}|=6$.\qed
 \end{lemma}

\begin{proof} Since $O^2(U)$
is generated by two conjugates of $Q_1J$, and $|J:C_J(Q_1)|= 3^2$ by
Lemma~\ref{CQ1}, we have that $|C_J(O^2(U))|\ge 3$. Since  the elements of order $5$ in $U$ act fixed-point-freely on $J_0$ we have  $C_J(O^2(U))\cap J_0 =1$. Thus $|C_J(O^2(U))|=3$ and, as $r_2$ centralizes
$J/J_0$ and normalizes $C_J(O^2(U))$, we get that $C_J(O^2(U))=
C_J(U)$. Since $|F:U|= 6$, $U$ is a maximal subgroup of $F$  and $t_1$ inverts $J$, we learn that $|C_J(U)^F|= |C_J(U)^{N_G(J)}|=6$.
\end{proof}}

\begin{lemma}\label{hyper1}
Suppose that $B \le J_0$ with $|B|=3^3$. Then $B$ contains a
conjugate of $Z$.
\end{lemma}

\begin{proof} Recall that $J_0$ is a non-degenerate quadratic space by Theorem~\ref{S6}(i). Hence this result follows because every subgroup of order $3^3$ in the $J_0$ contains a singular vector and the singular one-spaces in $J_0$ are $G$-conjugate to $Z$.
\end{proof}

We now fix some further notation.  First let $W =C_F(r_2)$. \blue{By the Frattini Argument we see that $WJ/J \cong C_{F/J}(r_2J)$ and  so}
$WJ/J \cong 2\times \Sym(4)$ and $J\cap W$ has index $3$ in $J$ by
Theorem~\ref{S6} (iv).

If $J = J_0$, set $\tau = 1$, whereas, if $J > J_0$,  select  $\tau
\in C_J(U)^\#$.

Suppose that $J > J_0$. Then $\t \not=1$. Let $$\mathcal T=
\t^F=\{\t_1=\tau, \dots, \tau_6\}$$ be the six $F$-conjugates of
$\t$. Then, as $[J,r_2]$ has order $3$ by Theorem~\ref{S6} (iv),
$r_2$ acts as a transposition on $\mathcal T$  and $r_2$ centralizes
$\tau$ (as $r_2 \in U$). Since $WJ/J \cong 2\times \Sym(4)$ and $W$
has orbits of length $2$ and $4$ on $\mathcal T$. It follows that,
after adjusting notation if necessary,
$\t^W=\{\t_1,\t_2,\t_3,\t_4\}$  and $\tau_5^{r_2}= \tau_6$. We
further fix notation so that $Q_1$ acts as $\langle
(\tau_2,\tau_3,\tau_4)\rangle$ and, since $r_1$ is conjugate to
$r_2$ in $N_G(J)$ and inverts $QJ/J$, we may suppose that $r_1$
induces the transposition $(\t_2,\t_3)$ on $\t^W$.

For $1\le i\le 4$, let $$J_i =\langle \t_j \mid 1 \le j\le 4, i\neq
j\rangle.$$ Then each $J_i$ is centralized by $r_2$  and is a
hyperplane of $C_J(r_2)$. Further  $$J_i \cap J_j = \langle \t_k
\mid 1\le k \le 4, k \not \in \{i,j\}\rangle.$$

Let $\rho \in [J,r_2]^\#$. Then $\rho \in (A\cap Q_2)\setminus Z$ as $|[J,r_2]| = 3$.
Since $[J,r_1]\le A \cap Q_1$,  we know  $[\rho,r_1]=1$. From the
choice of $\tau$ and $\rho$, we have that $\langle Q_1,r_1\rangle$
and $\langle \tau,\rho\rangle$ commute.

\blue{We now select and fix once and for all $$\rho \in [A \cap Q_2,r_2]^\#.$$}

For $J_0 = J$ we have to define the groups $J_1, J_2, J_3$ and $J_4$
differently. Set $J_1 = C_A(r_2)= A \cap Q_1$. So $J_1$ is
normalized by $\langle r_1,r_2, Q_1,J \rangle$ which has index $4$
in $W$. \blue{Observing that  $Z$ is centralized by the Sylow 3-subgroup $S$ of $F$ and $\langle W,S \rangle = F$,  yields that  $W$ is not contained in $M$.} As $Z$ is the unique
element of $Z^G$ contained in $J_1$, we have
$J_1^W=\{J_1,J_2,J_3,J_4\}$ and $W$ acts $2$-transitively on
$J_1^W$. As $r_1 \sim_M r_2$, all the elements in $J_1 \setminus Z$
are conjugate to $\rho$. Therefore,  as all the subgroups $J_i$ are
centralized by $r_2$, we have that $|J_i \cap J_j| = 3$ for all $i
\not= j$ and these intersections are conjugate to $\langle \rho
\rangle $. We capture some of the salient properties of these
subgroups in the next lemma.

\begin{lemma}\label{Jiprops} For $1\le i \le 4$, $J_i \le C_G(r_2)$ and $N_{N_G(J)}(J_i)$ contains a Sylow $3$-subgroup of $N_G(J)$.
\end{lemma}

\begin{proof}  If $J> J_0$, this is transparent from the construction of the subgroups. In the case that $J= J_0$, we have already mentioned that the subgroups commute with $r_2$.  Also we have $J_1 = A \cap Q_1$ is normalized by $S$ and as $J_i$, $2 \le i \le 4$ are conjugates of $J_1$ in $N_G(J)$, we have $N_{N_G(J)}(J_i)$ contains a Sylow $3$-subgroup of $N_G(J)$.
\end{proof}

Note also that when $|J|=3^5$,  $\rho \in \langle
[\t_5,r_2]\rangle$. It follows that $\langle \t_5,\t_6\rangle$
contains $\rho$ in this case. When $J=J_0$, of course we have
$\t_i=1$.  Thus to handle the two possible cases simultaneously we
will consider the group $\langle \t_5,\rho\rangle$.

\begin{lemma}\label{56cent}  $\langle \tau_5, \rho\rangle$ is centralized by $ JQ_1R_1$. Further $C_G(\langle \tau_5, \rho\rangle) \not \le M$.
\end{lemma}

\begin{proof} Set $X=\langle \tau_5,\rho\rangle$.  If $|J|= 3^4$, then $X=\langle \rho\rangle\le A \cap Q_2 $ and  the lemma holds. So suppose that $|J|=3^5$. Then $X = \langle\tau_5,\tau_6\rangle$ is centralized by $J$. Further,  as $\{\t_5,\t_6\}$ is a $W$-orbit and $Q_1 \le C_F(r_2) \le W$, $Q_1$ centralizes $X$.  Since $C_J(Q_1) = C_J(R_1)$ by Lemma \ref{CQ1} we now have $[X,R_1] = 1$ and this completes the proof.

Notice that $\langle \t_5, \rho\rangle$ is centralized by a subgroup
of index $2$ in $W$ and so $C_G(\langle \t_5,\rho\rangle)$ is not
contained in $M$.
\end{proof}

\begin{lemma}\label{CMtr}  The following hold.
\begin{enumerate}
\item $C_M(\rho) = JQ_1R_1\langle r_2t_1 \rangle$; and
\item If $J>J_0$, $C_M(\langle \tau_5,\rho \rangle)= JQ_1R_1$.
 \end{enumerate}
\end{lemma}

\begin{proof}  We calculate that  $C_M(\rho)$ contains $JQ_1R_1\langle r_2t_1\rangle$.  As $JQ_1R_1\langle r_2t_1\rangle$ covers $C_{M/Q}(\rho Z)$ (i) holds.

By Lemma~\ref{56cent},  $\langle \tau_5,\rho \rangle$ is centralized
by $JQ_1R_1$. Since, by Lemma \ref{NMJ}(iii), $r_2t_1$ conjugates $\langle \t_5 \rangle$ to $\langle \tau_6 \rangle$, part
(ii) follows from (i).
\end{proof}

\begin{lemma}\label{ZUN} $Z$ is the unique $G$-conjugate of $Z$ in $\langle \tau_5, \rho , Z\rangle$.
\end{lemma}

\begin{proof} Since $Z$ is weakly closed in $Q$, $Z$ is the unique conjugate of $Z$ in $\langle Z, \rho\rangle$. Also, as
$\t_5$ is not contained in $J_0$ and all the $G$-conjugates of $Z$
in $J$ are contained in $J_0$, there are no $G$-conjugates of $Z$ in
$\langle \t_5 , \rho , Z\rangle\setminus \langle\rho,Z\rangle$. This
proves the claim.
\end{proof}

\begin{lemma}\label{nots3}  Assume that $J>J_0$. Then $N_G(\langle r_1,r_2\rangle )/C_G(\langle r_1,r_2\rangle) \not \cong \Sym(3)$. \end{lemma}

\begin{proof} Let $U= \langle r_1,r_2\rangle$.   By a Frattini Argument $C_M(U)$ covers $M/Q$.  Hence as $J> J_0$, we have $|C_M(U)|_3=3^3$ and so
$D=C_J(U)$ is a Sylow $3$-subgroup of $C_M(U)$. Since $Z \le D$, we have $C_G(D) = C_M(D)= JU$ which is
$3$-closed. Therefore, $N_G(D) \le N_G(J)$. Since $r_1$ and $r_2$ act as transpositions on $\mathcal T$,
$|N_{F}(DU)/J|= 32$ and so we deduce that $D \in \Syl_3(C_G(U))$. Let $P = N_{N_G(U)}(D)$. Then by the
Frattini Argument $PC_G(U)= N_G(U)$. Therefore, if $N_G(U)/C_G(U) \cong \Sym(3)$, then $r_2$ and $r_1r_2$ are
conjugate in $P$. But $P \le N_G(J)$, $r_2 \in F\setminus F'$ and $r_1r_2\in F'$ which is a contradiction. Hence
$N_G(U)/C_G(U) \not \cong \Sym(3)$.
\end{proof}

\section{A further $3$-local subgroup and a $2$-local subgroup in the centralizer of an involution}

In this section we study the normalizer of $\langle\tau_5,
\rho\rangle$ and construct a $2$-local subgroup of $C_G(r_2)$.

\begin{lemma}\label{J_0signalizers}
We have $\signal_G(J_0,3')=\{1\}$.
\end{lemma}

\begin{proof}  Suppose that $1 \not= Y\in \signal_G(J_0,3')$. Then, as every hyperplane of $J_0$ contains a conjugate
of $Z$ by Lemma~\ref{hyper1}, \blue{and by coprime action $Y$ is generated by centralizers of hyperplanes  of $J_0$,} we may assume that $X= C_{Y}(Z)\neq
1$. So $X\in \signal_M(J_0,3')$. As $X$ is normalized by $A = J_0
\cap Q$ and $X$ normalizes $Q$, $$[A,X] \le Q \cap X=1.$$ \blue{But then $X$ centralizes a maximal abelian subgroup
of $Q$ and consequently $[Q,X]=1$ by Lemma~\ref{w=1}, which is a contradiction. Thus $\signal_G(J_0,3')=\{1\}$.}
\end{proof}

\begin{lemma}\label{notSp} Assume $J=J_0$. Then $C_G(\rho) \not\cong \langle \rho \rangle \times \Sp_6(2)$.
\end{lemma}

\begin{proof} Suppose that $C_G(\rho) \cong \langle \rho \rangle \times \Sp_6(2)$. Set $E = E(C_G(\rho))$. Then $E\cong \Sp_6(2)$. We have that $r_2$ inverts $\rho$ and  centralizes  $J/\langle \rho\rangle$, so as $J \cap E$ has order $3^3$ and $C_E(J\cap E)= J\cap E$, $r_2$ induces the trivial automorphism on $E$. Hence $N_G(\langle \rho \rangle) \cong \Sym(3) \times E$ and $[E,r_2] = 1$.
In $E \cap J$ there is an element  $\wt{\rho}$ with $N_E(\langle
\wt{\rho} \rangle) \cong \Sp_2(2) \times \Sp_4(2)$ (see \cite[page 46]{Atlas}). Hence
$N_{N_G(J)}(\langle \rho \rangle) \cap N_{N_G(J)}(\langle
\wt{\rho} \rangle)$ contains a Sylow 2-subgroup  $T$ of $N_G(J)$.
Now $\langle \rho, \wt{\rho} \rangle = C_J(i)$, where $i\in T'
\le F^\prime$.  \blue{Since $\mathrm O_2^-(3) \cong \Dih(8)$ and $N_G(J)/N_G(J)' \cong 2 \times 2$, we see that the involutions in $N_G(J)'/J$ invert a $-$-space and centralize a $+$-space with respect to the form given in Theorem~\ref{S6} (i).   In particular $C_J(i)$ is a $+$-space  and so $i$ centralizes a conjugate of }
$Z$.  Hence $\langle \rho, \wt{\rho}
\rangle$ contains a conjugate of $Z$. But $C_E(\langle \rho, \wt{\rho} \rangle)^{(\infty)} \cong \Sp_4(2)^{(\infty)} \cong \Alt(6)$ contradicts the fact that $M$ is soluble.
\end{proof}

\begin{lemma}\label {ctr} Let $B$ be a maximal subgroup of $\langle \tau_5, \rho, Z\rangle$ and  assume
that  $C_G(B) \not\leq M$.  Then $B\in \langle
\tau_5,\rho\rangle^{Q_2}$ and either
\begin{enumerate}
\item  $J> J_0$ and  $C_G(B) \cong B  \times \SU_4(2)$;  or
\item  $J=J_0$ and  $C_G(\rho) \cong \langle \rho \rangle \times \Aut(\SU_4(2))$.
\end{enumerate}
\end{lemma}

\begin{proof} Set $U=\langle Z, \t_5 ,\rho \rangle$, let $B$ be a maximal subgroup of  $U$,  $X = C_G(B)$  and  $\wt X = X/B$.
Assume that   $X \not \le M$.  By Lemma~\ref{ZUN}, $Z$ is the unique
conjugate of $Z$ in $U$ and so, as $C_G(B) \not \le M$, $U= ZB$
 and $\wt{N_X(Z)} = N_{\wt X} (\wt Z)$.

 Assume that $J> J_0$. Then, by Lemma \ref{CMtr}(ii), $N_X(Z) = X \cap M = JQ_1 R_1$ and so
 $N_{\wt X} (\wt Z) = \wt {N_X(Z)} = \wt {JR_1Q_1}\cong
3^{1+2}_+.\SL_2(3)$ which is isomorphic to the centralizer of a
$3$-central element in $\SU_4(2)$.  \blue{As $B \cap \langle \rho, Z \rangle \not= 1$, we may assume hat $\rho \in B$. Then by \ref{CMtr}(i) we have $C_M(B) \leq C_M(\rho) = JQ_1R_1\langle r_2t_1 \rangle$. As $|[U, r_2t_1]| = 9$, we get $C_M(B) \leq JQ_1R_1$ and so } $z \in Z^\#$ is not
$X$-conjugate to its inverse by Lemma \ref{CMtr},
$\signal_G(J_0,3')=\{1\}$ by Lemma~\ref{J_0signalizers} and $C_G(B)
\not \le M$, we may apply Hayden's Theorem \ref{Hayden} to get that
$\wt X \cong \SU_4(2)$.
 Finally, as $JQ_1$, splits over $B$, $X$ splits over $B$ by Gasch\"utz's Theorem \cite[9.26]{GLS2}. Hence  $X$ has the structure described in (i).

Assume that $J=J_0$. In this case $B$ is $Q_2$-conjugate to $\langle \rho \rangle$.
By Lemma \ref{CMtr}(i), $C_X(Z) = X \cap M = JQ_1 R_1$, \blue{as $r_2t_1$ inverts $Z$} and so $\wt
{C_X(Z)}$  is isomorphic to the centralizer of a $3$-central element
in $\SU_4(2)$. Since $r_2t_1$ inverts $z$,  we may use Prince's
Theorem \ref{PrinceThm} to obtain $\wt X \cong \Aut(\SU_4(2))$ or
$\Sp_6(2)$. Again Gasch\"utz's Theorem implies that $X \cong \langle
\rho \rangle \times E$ where $E \cong \wt{X}$.
Therefore, by Lemma \ref{notSp}, $X$ has the structure claimed  in
(ii).

Now we consider the possibilities for $B$ when $J>J_0$. We have  $B
\le U$ and $C_G(B) \not \le M$. Thus, by (i), $C_G(B)\cong B \times
E$ where $E \cong \SU_4(2)$. Consequently, $N_{C_G(B)}(J) \cong 3^2
\times (3^3{:}\Sym(4))$. Since $N_{C_G(B)}(J) \ge Q_1$ and since
there are exactly three subgroups isomorphic to $\Alt(4)$ which
contain a given $3$-cycle in $\Sym(6)$, we see that $B$ is
$Q_2$-conjugate to $\langle \tau_5,\rho\rangle$ as claimed.
\end{proof}

We now set $r=r_2$ and aim to determine $$K=C_G(r).$$

We will frequently use the following observation.

\begin{lemma}\label{5.4} $C_J(r)Q_1 $ is a Sylow $3$-subgroup of $K$.
\end{lemma}

\begin{proof} Certainly $C_J(r)Q_1\le K$ by Lemma~\ref{Qaction} (i). Because  $[Q_1,C_J(r),Q_1] = [A \cap Q_1, Q_1]= Z$, we have that $Z$ is a characteristic subgroup of $C_J(r)Q_1$ and so it follows that $N_K(C_J(r)Q_1) \le C_M(r)$. As $C_J(r)Q_1\in \syl_3(C_M(r))$, the lemma holds.
\end{proof}

Define $E= E(C_G(\langle \tau_5,\rho\rangle))$. Then  $E \cong
\SU_4(2)$ by Lemma~\ref{ctr}.

\begin{lemma}\label{normalizer} We have $E \langle t_1, \t_5\t_6\rangle\le K$ and $E \langle t_1 \rangle \cong \Aut(\SU_4(2))$.
\end{lemma}

\begin{proof} We know that $r$ inverts $\rho$ and exchanges $\t_5$
and $\tau_6$. Hence $r$ normalizes $B= \langle \tau_5,\rho\rangle$ and consequently $r$ normalizes $E$.
Furthermore, $r$ centralizes $J \cap E$ and since no involutory automorphism of $E$ acts in this way \blue{(see
\cite[page 26]{Atlas})}, we have that $r$ centralizes $E$. Therefore  $E \le K$.

Since $t_1$ inverts $J$,  $t_1$ normalizes $\langle \t_5,\rho\rangle$ and $t_1$ therefore normalizes $E$. Since
$t_1$ inverts $J \cap E$ \blue{and, by \cite[page 26]{Atlas}, no inner automorphism of $\SU_4(2)$ inverts an
elementary abelian group of order 27}, we have $E \langle t_1 \rangle \cong \Aut(\SU_4(2))$.
\end{proof}

From Lemmas~\ref{CMtr} and \ref{ctr}  we have $Q_1R_1 \le E$.
Furthermore, as $W(=C_F(r))$ normalizes $[J,r]= \langle
\rho\rangle$, we also have that $C_W(\rho) \le E$. In  particular,
we have

\begin{lemma}\label{Egen} $\langle \tau_5\tau_6\rangle E= \langle C_W(\rho), Q_1R_1C_J(r)\rangle $.
\end{lemma}
\begin{proof}  \blue{As $Q_1R_1C_J(r)$ contains the maximal parabolic subgroup of shape $3^{1+2}_+.\SL_2(3)$
of $\SU_4(2)$}, we have that $Y=Q_1R_1C_J(r)$ is a maximal subgroup of $E\langle \tau_5\tau_6\rangle$ and $C_W(\rho)\not\le Y$.

\end{proof}

When $J > J_0$, as $N_G(J)$ acts $2$-transitively on $\mathcal T$,
$\langle \t_5, \tau_6\rangle$ is $G$-conjugate to each subgroup $J_i
\cap J_j$ for $1\le i <j\le 4 $. When $J= J_0$ we have the same
result from the construction of $J_1, J_2, J_3$ and $J_4$ in Section
4. Hence we may apply Lemma~\ref{ctr} to obtain the following
conclusion.

\begin{lemma}\label{CJij} Assume that $1\le i < j \le 4$.\begin{enumerate}  \item If $J>J_0$, then $C_G(J_i\cap J_j) \cong
(J_i \cap J_j) \times \SU_4(2)$; and \item If $J=J_0$, then
$C_G(J_i\cap J_j) \cong (J_i \cap J_j) \times \Aut(\SU_4(2))$.
\end{enumerate}\qed
\end{lemma}

For $1\le i<j\le 4$, define  $$E_{ij}= E(C_G(J_i\cap J_j)).$$

\begin{lemma}\label{3central} For $1\le i < j \le 4$ and $k\in\{i,j\}$,  $E_{ij}\cap J_k$ is conjugate to $Z$ and is $3$-central in $E_{ij}$. In particular, $C_G(J_i) \cong  (J_i \cap J_j)\times 3^{1+2}_+{:}\SL_2(3)$ if $J>J_0$ and $C_G(J_i) \cong  (J_i \cap J_j)\times 3^{1+2}_+{:}\SL_2(3).2$ if $J =J_0$. \end{lemma}

\begin{proof} Let $1\le i\le 4$. Then by Lemma~\ref{Jiprops}, $J_i$ is normalized by a Sylow $3$-subgroup $T_i$ of $N_G(J)$ and $C_{T_i}(J_i)$ has
index $3$ in $T_i$. In particular, as $|C_{G}(J_i \cap J_j)|_3= 3|J|
$, we see that $C_{T_i}(J_i) \in \Syl_3(C_{G}(J_i \cap J_j))$.
Therefore $J_i \cap E_{ij}$ is normalized by a Sylow $3$-subgroup of
$E_{ij}$. As $|J_i\cap E_{ij}| =3$, we have that $J_i\cap E_{ij}$ is
$3$-central in $E_{ij}$ as $J_i$ is normal in $T_i$, we see that
this subgroup is also normal in a Sylow $3$-subgroup of $G$.
\end{proof}

 Define
$$ \Sigma=\langle O_2(C_K(J_k))\mid 1\le k\le 4\rangle.$$

In the next lemma we use the fact that if $x \in \SU_4(2) = X$ is an
involution which centralizes a subgroup of order $9$, then $x$ is
$2$-central  and $$C_X(x) \cong 2^{1+4}_+.(3\times \Sym(3))\cong
(\SL_2(3)\circ \SL_2(3)).2$$ where $\circ$ denotes a central product
(see \cite[page 26]{Atlas}).

\begin{lemma}\label{extraspec} Assume that $1 \le i < j \le 4$. Then
\begin{enumerate}
\item $O_2(C_K(J_i)) \cong O_2(C_K(J_j)) \cong \Q_8$, $ [O_2(C_K(J_i)),O_2(C_K(J_j)) ]=1$ and  $O_2(C_K(J_i\cap J_j))= O_2(C_K(J_i)) O_2(C_K(J_j))\cong 2^{1+4}_+;$ and
\item $\Sigma$ is extraspecial of $+$-type  and order $2^9$.
\end{enumerate}
\end{lemma}
\begin{proof}
Suppose that $1\le i< j \le 4$. Then $J_i \le C_G(r)$ by Lemma~\ref{CJij}.  If $J> J_0$, we have $r \in E_{ij}$
by Lemma~\ref{ctr}. If $J= J_0$, then \blue{as $[J_1,R_2] = 1$ and $r \in Z(R_2)\le C_G(J_1)^\prime$ we have} $r
\in E_{12}$ and consequently $r \in E_{ij}$ as $W$ acts $2$-transitively on $\{J_1,J_2,J_3, J_4\}$.

Since  $r\in E_{ij}$ and $|C_J(r)\cap
 E_{ij} |_3\ge 9$,  $r$ is a $2$-central involution in
$E_{ij}$. It follows that $K \cap E_{ij} $ has shape
$2^{1+4}_+.(3\times \Sym(3))$ and, in particular, $O_2(C_K(J_i\cap
J_j))\cong 2^{1+4}_+.$ Furthermore, as $J_i \cap E_{ij}$ is
$3$-central by Lemma \ref{3central},  we get $
O_2(C_K(J_i))\cong \Q_8$ and $O_2(C_K(J_i \cap J_j)) =
O_2(C_K(J_i))O_2(C_K(J_j))$. Since $O_2(C_K(J_i \cap J_j))$ contains
exactly two subgroups isomorphic to $\Q_8$, we have that
$[O_2(C_K(J_i)),O_2(C_K(J_j))]=1$. This completes the proof of  (i).

Part (i) shows that  $\Sigma$  is isomorphic to a central product of
$4$ quaternion groups. Hence $\Sigma$ is extraspecial of $(+)$-type
and  order $2^9$. So (ii) holds.
\end{proof}

Recall from Corollary \ref{cor} and \ref{NMSagain}, $t_2 \in N_G(S)
\le M \cap N_G(J)$ and $R_1^{t_2}= R_2$.

\begin{lemma}\label{R2} We have $J_1 $ is centralized by $R_2$, $R_2 \leq \Sigma$ and $R_2 = C_{\Sigma}(Z)$.
\end{lemma}

\begin{proof}  Suppose first that $J=J_0$.  Then  $J_1 = C_A(r) \le Q_1 = C_Q(R_2)$  by Lemma~\ref{Qaction} (i).  So $[J_1,R_2] = 1$. Hence
$R_2 = O_2(C_K(J_1)) \leq \Sigma$.

Assume that $J>J_0$.  We have that $\tau_1$ commutes with $Q_1$ and
$[\langle \t_5, \t_6 \rangle,Q_1] = 1$ by Lemma \ref{CMtr}. Hence
$C_J(Q_1)= \langle \t_1,\t_5,\t_6\rangle= \langle \tau_5, A\cap
Q_2\rangle$. Thus $C_{J}(Q_2)= C_J(Q_1)^{t_2}= \langle
\t_2,\t_3,\t_4\rangle= \langle\t_2,A\cap Q_1\rangle$. By
Lemma~\ref{CQ1}, $C_J(Q_1)$ is centralized by $R_1$, thus $J_
1=\langle \t_2,\t_3,\t_4\rangle$ is centralized by $R_2=R_1^{t_2}$.
Hence $R_2 = O_2(C_K(J_1)) \leq \Sigma$.

Since $R_2$ commutes with $Z$, we have $R_2 \le C_{\Sigma}(Z)$ and, as $C_{\Sigma}(Z)$ is extraspecial we have
that $R_2 = C_{\Sigma}(Z)$ from the structure of $M$.
\end{proof}

\begin{lemma}\label{WinK}  We have $W\langle t_1\rangle \le N_K(\Sigma)$.
\end{lemma}
\begin{proof}Since $W\langle t_1 \rangle$ permutes $ \{J_1, J_2,J_3,J_4\}$  and is contained in $K$, $W\langle t_1
\rangle \le N_K(\Sigma)$ by the definition of $\Sigma$.\end{proof}

\begin{lemma}\label{NCJr} We have $W = N_{{K}}(C_J(r)) = N_{N_K(\Sigma)}(C_J(r))$. In particular $N_{N_K(\Sigma)}(C_J(r))$ controls $K$-fusion in $C_J(r)$.
\end{lemma}

\begin{proof} We have that $C_G(C_{J}(r)) = J\langle r \rangle$. Hence $J$ is normal in $N_G(C_{J}(r))$. Now we have that $W = N_K(C_{J}(r))$. By Lemma
\ref{WinK} we have $W \leq N_K(\Sigma)$ and so $N_K(C_{J}(r)) =
N_{N_K(\Sigma)}(C_{J}(r))$. Further by Lemma \ref{J ab} we have that
$N_G(J)$ controls fusion in $J$ and so $N_K(C_{J}(r))$ controls
$K$-fusion in $C_{J}(r)$. As $N_K(C_{J}(r)) = N_{N_K(\Sigma)}(C_{J}(r))$ this
fusion takes place in $N_K(\Sigma)$.
\end{proof}

\begin{lemma}\label{Jsig}  Every  $J_1$-signalizer in $K$ is contained in $\Sigma$. In particular, $N_K(J_1) \le N_K(\Sigma)$.
\end{lemma}

\begin{proof}
Let $\Sigma_1 \leq K$ be a $J_1$-signalizer. Let $X_1$ be a hyperplane in $J_1$ such that $C_G(X_1) \leq M$. Then
$C_{\Sigma_1}(X_1) \leq M$ is normalized by $J_1$ and so $C_{\Sigma_1}(X_1)\le R_2 \le \Sigma$ by Lemma \ref{R2}. In
particular $[C_{\Sigma_1}(Z), J_1] = 1$.

Suppose next that $X_1$ is a hyperplane such that $C_G(X_1) \not\leq M$. Then, by Lemma \ref{ctr}, we may assume
that $X_1= J_1 \cap J_2$.  Since $r$ is $2$-central in $E_{12}$, $O_2({C_K(J_1 \cap J_2)})$
 is the  unique maximal $J_1$-signalizer in
$C_G(X_1)$. Hence by Lemma~\ref{extraspec} (i) we have that
$C_{\Sigma_1}(X_1) \leq \Sigma$ in this case as well. Because
$$\Sigma_1 = \langle C_{\Sigma_1}(X_1)\mid |J_1 : X_1| \leq 3
\rangle \leq \Sigma,$$ we have that every $J_1$-signalizer is contained in $\Sigma$. Thus $\Sigma$ is the unique
maximal member of $\signal_K(J_1,3')$ and so $N_K(J_1)\le N_K(\Sigma)$ \blue{as $N_K(J_1)$ acts via conjugation
on the maximal elements of $\signal_K(J_1,3^\prime)$.}
\end{proof}

\begin{lemma}\label{faith} $C_K(\Sigma)= \langle r\rangle$.
\end{lemma}

\begin{proof} If $C_K(\Sigma)$ is a $3'$-group,
then $C_K(\Sigma)$ is normalized by $J_1$ and so $C_K(\Sigma) \le Z(\Sigma)=\langle r\rangle$ by
Lemma~\ref{Jsig}. So suppose that $C_K(\Sigma) $ has order divisible by $3$. Since by Lemma \ref{5.4}  $C_J(r)Q_1 \in \Syl_3(K)$ and
$C_J(r)Q_1\le W \le N_K(\Sigma) $  by Lemma~\ref{WinK}, we have $C_J(r)Q_1\cap C_G(\Sigma)$ is a Sylow
$3$-subgroup of $C_G(\Sigma)$. As $Z $ does not centralize $\Sigma$, we have $ C_J(r)Q_1\cap C_G(\Sigma) \le
C_J(r)$. Now, for $1\le i<j\le 4$
$$C_{C_J(r)}(O_2(C_K(J_i \cap J_j)) )= J_i \cap J_j$$ and
consequently $C_{C_J(r)}(\Sigma) \le J_1\cap J_2 \cap J_3 \cap J_4=1$ which is a contradiction.
\end{proof}

\begin{lemma}\label{irred} $\Sigma/\langle r \rangle$ is a minimal normal subgroup of $N_K(\Sigma)/\langle r\rangle$.
\end{lemma}

\begin{proof} Suppose that $U\le \Sigma$ and $U/\langle r \rangle$ is a minimal normal subgroup of $N_K(\Sigma)/\langle r\rangle$ of minimal order. Aiming for a contradiction, assume that $U \neq \Sigma$.
 Then either $|\Sigma: U| \le 2^4$ or $|U/\langle r \rangle|\le 2^4$. In particular, as $Q_1$ normalizes $\Sigma$ (see Lemma \ref{Jsig}) and $\GL_4(2)$ has elementary abelian Sylow $3$-subgroups,   $Z$ centralizes one of  $U$ or $\Sigma/U$. By Lemma~\ref{R2}, either $U \le R_2$ or $|\Sigma:U| \le 2^2$ and $U \ge [\Sigma,Z]$.

  Since $C_{J}(r)$ acts non-trivially on $R_2$, we get $U= R_2$ or $U= [\Sigma,Z]$. In the latter case, we have $U_1= C_{\Sigma}(U)$ is normalized by $N_K(\Sigma)$  and has order smaller than $U$. Hence the minimal choice of $U$ implies that $U=R_2$. However $W \le N_G(\Sigma)$ by Lemma~\ref{WinK} and $W$ does not normalize $R_2$ and so we have a contradiction.
\end{proof}

\begin{theorem}\label{Normalizer} One of the following holds.
\begin{enumerate} \item $J=J_0$ and $N_G(\Sigma)/\Sigma \cong \Aut(\SU_4(2))$ or $\Sp_6(2)$; or
\item $J > J_0$ and  $N_G(\Sigma)/\Sigma \cong (3 \times \SU_4(2)){:}2$.
\end{enumerate} Furthermore, $E \langle \t_5\t_6,t_1\rangle \le N_K(\Sigma)$ and  $\Sigma/\langle r\rangle$ is isomorphic to the natural $E\Sigma/\Sigma$-module.
\end{theorem}

\begin{proof} From Lemma \ref{WinK} we have that $W\langle t_1\rangle \leq N_G(\Sigma)$. Set $L = J_1Q_1$. Then
$L \leq W$ and so $L \leq N_G(\Sigma)$. By Lemma \ref{Jsig} we have that $\Sigma$ is a maximal signalizer in $K$
for $L$ and for $C_J(r)$. Hence $N_K(L)$  and $N_K(C_{J}(r)) $ both normalize  $\Sigma$.

Suppose that $J= J_0$.Then $J_1Q_1 = (A\cap Q_1)Q_1\le Q_1$ and so $R_1\le N_K(Q_1) \le N_K(\Sigma)$. Therefore Lemma~\ref{Egen} implies that $\langle E, t_1\rangle \le N_K(\Sigma)$.  In particular, we have $C_{N_K(\Sigma)/\Sigma}(Z\Sigma /\Sigma)$ is isomorphic to the centralizer of a $3$ element in $\SU_4(2)$ and is inverted by $t_1\Sigma$. Hence Theorem~\ref{PrinceThm} shows that (i) holds.

Suppose that $J > J_0$. This time $N_K(J_1Q_1)$ does not contain $R_1$.  On the other hand $N_{K}(\Sigma) \ge N_K(C_J(r)) \Sigma=  W\Sigma$ and $W\Sigma/\Sigma $ has shape $3^4{:}(\Sym(4) \times 2)$.  By the Frattini Argument, $N_{N_K(\Sigma)/\Sigma}(C_J(r)\Sigma/\Sigma)= N_{N_K(\Sigma)}(C_J(r))$. Since $N_K(C_J(r)) = W$, we now have $N_{N_K(\Sigma)/\Sigma}(C_J(r)\Sigma/\Sigma) = W\Sigma/\Sigma$.

Since $C_G(\Sigma) = \langle r \rangle$ by Lemma~\ref{faith}, we have that $N_K(\Sigma)/\Sigma$ is isomorphic to
a subgroup of $\mathrm O_8^+(2)$.  Because $N_{N_K(\Sigma)/\Sigma}(C_J(r)\Sigma/\Sigma) = W\Sigma/\Sigma$, we infer from the list of maximal subgroups of $\mathrm O_8^+(2)$ given in \cite[page 85]{Atlas} that either  $N_K(\Sigma)= W\Sigma$ or $N_K(\Sigma)/\Sigma \cong (3\times \SU_4(2)){:}2$. In the latter case we have (ii) so suppose that $N_K(\Sigma)= W\Sigma$.  Let $T \in \Syl_2(N_K(\Sigma))$. We claim that $T \in \syl_2(K)$. Assume that $x \in N_K(T) \setminus N_K(\Sigma)$. Then, as $\Sigma^x \not=\Sigma$,   $J(T/\langle r \rangle) \not \le \Sigma/\langle r\rangle$. Hence, setting  $H= \langle J(T)^{N_K(\Sigma)}\rangle$  and noting that $|O_3(N_K(\Sigma)/\Sigma)|=3^4$, we may apply \cite[(32.5)]{Aschbacher} to get that $H/\Sigma$ is a direct product of four subgroups isomorphic to $\SL_2(2)$. But then the $2$-rank of $W/\Sigma$ is at least $4$ contrary to $T/\Sigma \cong \Dih(8) \times 2$. Hence $N_K(T) \le N_K(\Sigma)$ and, in particular, $T \in \Syl_2(K)$.

From Lemma~\ref{normalizer}, we have $E \le K$. Since $T\in \Syl_2(K)$, $T/\Sigma \cong \Dih(8)\times 2$ and $E$ contains an extraspecial subgroup of order $2^5$ with centre $\langle r_1\rangle$, we have that $r_1$ is $K$-conjugate to an element of $\Sigma$. Thus there is some $x \in K$ such that $\langle r_1,r \rangle \le \Sigma^x$. Since $r_1^{t_2}= r$ and since  $r_1$ and $rr_1$  are  $\Sigma^x$-conjugate, we have $N_G(\langle r_1,r\rangle)/C_G(\langle r_1,r\rangle) \cong \Sym(3)$. This contradicts Lemma~\ref{nots3}. Hence (ii) holds.

We have already seen that $E \le N_K(\Sigma)$ if $J= J_0$. If $J>J_0$, then we have $N_{N_K(\Sigma)}(Z) $ contains a subgroup $(3 \times 3^{1+2}_+).\SL_2(3).2$. Since $N_K(Z) = C_M(r) = Q_1R_1R_2C_J(r)\langle t_1\rangle$, we have $C_M(r) \le N_K(\Sigma)$. Now $E \langle \t_5\t_6,t_1\rangle \le N_K(\Sigma)$ by Lemma~\ref{Egen}. Finally, as $E$ acts irreducibly on $\Sigma/\langle r\rangle$ by Lemma~\ref{irred}, we have that $\Sigma/\langle r\rangle$ is the natural $E$-module. \end{proof}

We need just  two final details before we can move on to determine the structure of $K$.

\begin{lemma}\label{containinK} The following hold.
\begin{enumerate}
\item $N_K(Z) \le N_K( \Sigma)$; and
\item $N_K(J_i \cap J_j) \leq N_K(\Sigma)$, for $1 \le i < j \le 4$.
\end{enumerate}
\end{lemma}

\begin{proof}  For (i) we note that $N_K(Z)= C_M(r) \le E \langle \t_5\t_6, t_1\rangle \Sigma\le N_K(\Sigma)$ by Theorem~\ref{Normalizer}.

By Lemma \ref{extraspec} (i) we have that
$O_2(C_K(J_i \cap J_j)) \leq \Sigma$ and, as $r$ is a $2$-central element in $E_{ij}$, $C_J(r) \in \Syl_3(C_K(J_i\cap J_j))$. Hence$$N_K(J_i \cap J_j) = N_{N_K(J_i \cap J_j)} (C_J(r)) O_2(C_K(J_i \cap J_j)) \le N_K(\Sigma)$$ by Lemma~\ref{Jsig}.
\end{proof}

\section{The structure of $K$}

In this section we prove Theorem~\ref{structureK} which asserts that
$K = N_K(\Sigma)$. We continue the notation introduced in the
previous sections. We further set $K_1 = N_K(\Sigma)$ and
denote by $\wt{\;}$  the natural homomorphism from $K$ onto $K/\langle r \rangle$.\\

By Lemma~\ref{irred}, the subgroup  $\wt \Sigma$  can be regarded as the $8$-dimensional
irreducible $\GF(2)$-module for $\wt K_1/\wt\Sigma$. Thus we may
employ the results of Proposition~\ref{fact} to obtain information
about various centralizers of elements of order $2$ and $3$ in $\wt \Sigma$.
Using  Proposition~\ref{fact}(ii), we have $\wt K_1$ has two orbits
on $\wt \Sigma$.
 We pick representatives   $\wt{x}$ and $\wt{y}$ of these orbits with $\wt x$ singular and $\wt y$ non-singular.  It follows that $x$ is an involution and $y$ has order 4.

Our aim is to show that $\wt \Sigma$ is strongly closed in $\wt K $
and then use Goldschmidt's Theorem \cite{Goldschmidt} to show that
$K= K_1$. We now begin the proof of Theorem~\ref{structureK}.

\begin{lemma}\label{y} We have $\wt{K}_1$ contains a Sylow 2-subgroup of $C_{\wt{K}}(\wt{y})$. In particular
$|C_{\wt{K}}(\wt{y})|_2 = 2^{12}$ if $E(\wt K_1/\wt \Sigma) \cong
\SU_4(2)$  and $|C_{\wt{K}}(\wt{y})|_2 = 2^{14}$ if
$\wt{K}_1/\wt{\Sigma} \cong \Sp_6(2)$.
\end{lemma}

\begin{proof}  Let $T$ be a Sylow 2-subgroup of $C_{\wt{K}_1}(\wt{y})$ and assume that $T_1$ is a $2$-group with $|T_1 : T| = 2$. Choose $u \in T_1 \setminus T$.
 If  $|\wt{\Sigma}^{u}\wt{\Sigma}/\wt{\Sigma}| \leq 2$, then $|\wt{\Sigma}^{u} \cap \wt{\Sigma}| \geq 2^7$.
But by Proposition~\ref{fact} (iv),  $\wt{K_1}$ has no $2$-elements
not in $\wt \Sigma$ which centralize a subgroup of index two in
$\wt{\Sigma}$. Therefore $\wt \Sigma= \wt \Sigma^u$ and so $u \in
T_1\cap \wt K_1= T$ which is a contradiction. Hence
$|\wt{\Sigma}^{u}\wt{\Sigma}/\wt{\Sigma}| \geq 4$.

If $E(\wt{K}_1/\wt{\Sigma} ) \cong \SU_4(2)$, then  $T/\wt{\Sigma}$
is a semidihedral group of order 16 by Proposition~\ref{fact}(ii).
Since $\wt{\Sigma}^{u}\wt{\Sigma}/\wt{\Sigma}$   is a normal
elementary abelian subgroup {of $T/\wt{\Sigma}$} of order at least
$4$, we have a contradiction. Hence $\wt{K}_1/\wt{\Sigma} \cong
\Sp_6(2)$ by Lemma \ref{Normalizer}.  Now Proposition~\ref{fact}(ii), gives
$$C_{\wt{K}_1}(\wt{y})/\wt{\Sigma} \cong \G_2(2).$$
 Since, by \cite[Table 3.3.1]{GLS3}, $\G_2(2)$ does not contain elementary abelian subgroups of order $16$,
$2^6\geq |\wt{\Sigma}^{u} \cap \wt{\Sigma}| \geq 2^5$. But then all
involutions in $\wt{\Sigma}^{u}$  centralize a subgroup of order at
least  $2^5$ in $\wt{\Sigma}$, and so Proposition~\ref{fact} (i) and
(iv) shows that all the involutions in
 $\wt{\Sigma}^{u}\wt{\Sigma}/\wt{\Sigma}$ are unitary transvections and are conjugate in $\wt{K_1}/\wt \Sigma$.  Since
the two classes of involutions in
$C_{\wt{K}_1}(\wt{y})/\wt{\Sigma} \cong \G_2(2)$ are not fused in
$\wt{K_1}/\wt{\Sigma}$, we infer that
$$\wt{\Sigma}^{u}\wt{\Sigma}/\wt{\Sigma} \leq
(C_{\wt{K}_1}(\wt{y})/\wt{\Sigma})^\prime\cong \G_2(2)'\cong
\SU_3(3).$$
 Since, by \cite[Table 3.3.1]{GLS3},  $\SU_3(3)$ has no elementary abelian groups of order $8$, we have $|\wt{\Sigma}^{u}\wt{\Sigma}/\wt{\Sigma} |=4$. This means that $|\wt{\Sigma}^{u} \cap \wt{\Sigma}| =2^6$ and consequently all the involutions in $\wt{\Sigma}^{u}\wt{\Sigma}/\wt{\Sigma} $ have the same centralizer in $\wt{\Sigma}$.
As centralizers of involutions in $\G_2(2)^\prime$ are maximal
subgroups \cite[page 14]{Atlas}, we conclude that $ \wt{\Sigma}^{u}
\cap \wt{\Sigma}$ is normalized by
$(C_{\wt{K}_1}(\wt{y})/\wt{\Sigma})^\prime$. Thus
$(C_{\wt{K}_1}(\wt{y})/\wt{\Sigma})^\prime$ centralizes $\wt
\Sigma$ which is impossible. This contradiction proves the lemma.
The order of $T$ is calculated from Proposition~\ref{fact}(iii).
\end{proof}

\begin{lemma}\label{3closed} Let $S_1$ be a Sylow 3-subgroup of $C_{\wt{K}_1}(\wt{x})$ or
$C_{\wt{K}_1}(\wt{y})$. Then $N_{\wt{K}}(S_1) \leq \wt{K}_1$.
 In particular, for $z \in \wt \Sigma^\#$, $C_{\wt{K_1}}(z)$ contains a
 Sylow $3$-subgroup of $C_{\wt K}(z)$.
\end{lemma}

\begin{proof} We consider  $\wt{y}$ first. By Proposition~\ref{fact}(iii), $S_1$ has centre of order $3$  and, as faithful $\GF(2)$-representations of extraspecial groups of type $3^{1+2}_+$  have dimension $6$, we have $|C_{\wt \Sigma}(Z(S_1))|=4$. As $|C_{\wt{\Sigma}}(Z)| = 4$ we may assume using Proposition \ref{fact}(iii) that $Z = Z(S_1)$. On the other,  Lemma~\ref{containinK} (i)  gives $C_M(r) \leq K_1$. Hence we have that
$N_{\wt{K}}(S_1) \leq \wt{K}_1$.

Now we consider $\wt{x}$. By Lemma \ref{extraspec} (i), we have
$O_2(C_K(J_1 \cap J_2)) \leq \Sigma$. Hence, comparing orders, we may assume that $S_1
= J_1 \cap J_2$. But then by Lemma~\ref{containinK} (ii)
$N_{\wt{K}}(S_1) \leq \wt{K}_1$.
\end{proof}

Let $\wt{E} \leq \wt{K}_1$ such that $\wt{E}/\wt{\Sigma} =
E(\wt{K}_1/\wt{\Sigma})$. We have that $\wt{E}/\wt{\Sigma}
\cong \SU_4(2)$ or $\Sp_6(2)$. By Proposition~\ref{fact}(iii)
there are exactly three classes of elements of order three in
$\wt{E}$. As a Sylow 3-subgroup of $\wt{E}$ is isomorphic
to $3 \wr 3$, there is a unique elementary abelian subgroup of
order 27, and this subgroup contains elements from each of the conjugacy classes of elements of order $3$. As
$\wt {C_{J}(r)} \cap \wt{E}$ is elementary abelian of order 27, there
are representatives of these elements in $\wt{C_{J}(r)} \cap
\wt{E}$. It follows that every element of order $3$ in $\wt K$ is
conjugate to an element of $\wt {C_J(r)}$. So using Lemma \ref{NCJr}
get the following lemma.

\blue{
\begin{lemma}\label{conj} Two elements of order three in $\wt{K}_1$ are conjugate in $\wt{K}$ if and only
 if  they are conjugate in $\wt{K}_1$.\qed
\end{lemma}

Supposing that $J > J_0$, we establish some further notation. Let $\sigma \in \wt {K_1}$  have order $3$ and $\sigma \wt \Sigma$  be  centralized by $\wt E/\wt \Sigma$. Then  $\sigma$ is not $\wt K$-conjugate to
any element in $\wt{E}$ by Lemma~\ref{conj}.
}

\begin{lemma}\label{3conj} Suppose that  $\wt{u} \in \wt{K}_1 \setminus \wt {\Sigma}$ is  an involution which is $\wt {K}$ -conjugate to some involution in
$\wt {\Sigma}$. Assume that  $\nu \in C_{\wt {K_1}}(\wt{u})$ is
an element of order three.  Then we have
\begin{itemize}
\item[(i)] $C_{\wt {\Sigma}}(\nu) \not= 1$;
\item[(ii)]   $\langle \nu \rangle \not\sim Z$ in $\wt {K}$;
\item[(iii)]    if $J=J_0$, then $\nu \not\sim \rho$ in $\wt  K$;
and
\item[(iv)] $|C_{\wt {E}}(\wt{u})|$ is not divisible by $9$.
\end{itemize}

\end{lemma}
\begin{proof} Let $\wt{a} \in \wt {\Sigma}$ with $\wt{a} \sim_{\wt K}\wt {u}$. By Lemma \ref{3closed}, $\wt{K}_1$ contains a Sylow 3-subgroup of
$C_{\wt {K}}(\wt{a})$. By Lemma \ref{conj}, $\nu$ is conjugate to
an element $\mu$ of $C_{\wt{K_1}}(\wt{a})$ inside of $\wt{K}_1$.
Now obviously $C_{\wt {\Sigma}}(\mu) \not= 1$ and so the same holds
for $\nu$ which is (i).
\\
If $\langle \nu \rangle $ is conjugate to $Z$ in $\wt {K}$ or to
$\langle \rho \rangle$ in case of $\tau = 1$,   this happens also in
$\wt{K}_1$ by Lemma \ref{conj}. Hence we may assume that $\wt{a}
$ is conjugate to $ \wt{u}$ in $M \cap K$, or $N_K(\langle \rho
\rangle)$,  which both are contained in $K_1$ by Lemma~\ref{containinK}, a contradiction. Hence also (ii) and (iii) hold.
\\
Assume now that $S_1 \leq C_{\wt {E}}(\wt{u})$, $|S_1| = 9$. Then
$S_1$ is conjugate into a Sylow 3-subgroup  $S_2$ of $C_{\wt
{E}}(\wt{a})$. So by Lemma \ref{3closed} and
Proposition~\ref{fact}(ii) we may assume that $\wt{a} =
\wt{y}$ and thus $S_2$ is extraspecial of order 27. Hence $S_1$
contains some element which is conjugate into $Z(S_2)$. But $Z(S_2)$
is conjugate to $Z$, and this contradicts (ii). This finishes the
proof.
\end{proof}

\begin{lemma}\label{2conj} Suppose that  $\wt{u} \in \wt{K}_1 \setminus \wt{\Sigma}$ is an involution which is $\wt{K}$-conjugate   to some involution in $\wt{\Sigma}$. Then either
\begin{enumerate}
\item $\wt{u} \in \wt{E}$, $|[\wt{\Sigma}, \wt{u}]| = 4$ and   $C_{\wt{E}}(\wt{u})$ has order $2^{13}$  if  $E(\wt K_1/\wt \Sigma ) \cong \SU_4(2)$  and  order $2^{15}$ if $\wt K_1/\wt \Sigma \cong \Sp_6(2)$; or
 \item  $J>J_0$,  $\sigma^{\wt{u}} = \sigma^{-1}$ and
 $C_{\wt{E}/\wt{\Sigma}}(\wt{u}) \cong 2 \times \Sym(4) \le \Sym(6)$, and $|[\wt{\Sigma}, \wt{u}]| = 16$.
     \end{enumerate}
\end{lemma}
\begin{proof} If $|[\wt{u},\wt{\Sigma}] = 16$, then all
involutions in $\wt{\Sigma}\wt{u}$ are conjugate by elements
of $\wt{\Sigma}$.  Hence, by Proposition~\ref{fact}(i),
$\wt{u}$ centralizes some non-trivial 3-element  $\nu  \in \wt{E}$. By
Lemma \ref{3conj}(i), $C_{\wt{\Sigma}}(\nu)  \not= 1$. If $J=J_0$, then by Proposition~\ref{fact}(iii) $\langle \nu \rangle$ is
conjugate to $Z$ or  $\langle \rho \rangle$, which contradicts Lemma
\ref{3conj} (ii),(iii). So assume  that $J> J_0$. If
$\wt{u} \not\in \wt{E}$, we have the assertion (ii) with
Proposition~\ref{fact}(i) and Lemma \ref{3conj}(iv).  So  assume
$\wt{u} \in \wt{E}$. Then
$C_{\wt{E}/\wt{\Sigma}}(\wt{u})$ is contained in a
parabolic subgroup of $\wt E /\wt \Sigma$  of shape $2^4{:}\Alt(5)$
and so $\nu$ acts fixed point freely on $\wt{\Sigma}$,
contradicting Lemma~\ref{3conj} (i).
\\
So assume that  $|[\wt{u},\wt{\Sigma}]| = 4$.  Then, by
Proposition \ref{fact} (v),
$C_{\wt{E}/\wt{\Sigma}}(\wt{u}\wt{\Sigma})$ has orbits of
length 1,6 and 9 on
$C_{\wt{\Sigma}}(\wt{u})/[\wt{\Sigma},\wt{u}]$. Hence there
are  exactly three conjugacy classes of involutions in
$\wt{\Sigma}\wt{u}$,   two of which have representatives
centralized by an element of order three. Assume that $\wt{u}$ is
one of these. Let $\hat{u}$ be the involution, which is  centralized
by $S_1$, a preimage of a Sylow 3-subgroup  of
$C_{\wt{K}_1/\wt{\Sigma}}(\wt{u}\wt{\Sigma})$. Set $S_2 =
C_{S_1}(\wt{u})$.
 Then, using  Lemmas \ref{conj},  \ref{3closed}
and  \ref{3conj}(iv), we see that $|S_2| = 3$. Therefore $\wt{u}
\not\sim \hat{u}$. In particular we have  $\wt{u} =
\hat{u}\wt{s}$ where  $\wt{s} \in C_{\wt{\Sigma}}(\wt{u})
\setminus [\wt{\Sigma},\wt{u}]$. Hence
$C_{C_{\wt{\Sigma}}(\wt{u})/[\wt{\Sigma},\wt{u}]}(S_2) \not=
1$. By Proposition \ref{fact} (vi), we get  $|C_{\wt{\Sigma}}(S_2)|
= 4$. So by Proposition \ref{fact}(iii) $S_2$ does not centralize involutions in $\Sigma$. Thus we
may assume that $S_2 = Z$. But this contradicts Lemma \ref{3conj}
(iii) and proves the lemma.
\end{proof}

\begin{lemma}\label{weakcl} We have $\wt{y}^{\wt{K}} \cap \wt E \subseteq \wt{\Sigma}$.
\end{lemma}

\begin{proof}
Assume $\wt{y} \sim_{\wt K} \wt{u}$ for some involution
$\wt{u} \in \wt E\setminus \wt{\Sigma}$. By Lemma \ref{y}, we
have that $|C_{\wt{K}}(\wt{y})|_2 = 2^{12}$ if $\wt E/\wt \Sigma
\cong \SU_4(2)$ or $2^{14}$ if $\wt E/\wt\Sigma \cong \Sp_6(2)$.
This conflicts with that information given in  Lemma \ref{2conj}.
Hence no such elements exist.
\end{proof}

\begin{lemma}\label{K1Syl2}  $\wt \Sigma$ is weakly closed in $\wt K_1$. In particular,
$\wt K_1$ contains a Sylow $2$-subgroup of $\wt K$.
\end{lemma}

\begin{proof} Assume that $T \in \Syl_2(\wt K_1)$, $w \in \wt K$ and  $\wt \Sigma^w \le T$ with  $\wt \Sigma \not= \wt \Sigma^w$. Then $\wt \Sigma^w\cap \wt E$ has order at least $2^7$ and therefore is generated by conjugates of $\wt y$. Thus Lemma~\ref{weakcl} implies that $\wt \Sigma^w\cap \wt E\leq \wt \Sigma$. But then $|\wt\Sigma^w\wt \Sigma/\wt\Sigma|=2$  and $|\wt \Sigma \cap \wt  \Sigma^w|=2^7$. Since $\wt K_1$ does not contain transvections, we have a contradiction.
\end{proof}

\begin{lemma}\label{tildex} No element of $\wt \Sigma$ is $\wt
K$-conjugate to an involution  $\wt{u} \in \wt{K}_1$ with
$|[\wt{\Sigma}, \wt{u}]| = 4$.
\end{lemma}

\begin{proof} Assume the statement is false.
Then, by Lemma \ref{weakcl} $\wt{u} \sim_{\wt K} \wt{x}$. Let
$T_1$ be a Sylow 2-subgroup of $C_{\wt{K}_1}(\wt{u})$ and $T_2$
be a Sylow 2-subgroup of $C_{\wt{K}}(\wt{u})$ with $T_1 \leq T_2$. By
Lemma~\ref{fact}(ii), \ref{2conj} and \ref{K1Syl2}, $|T_2 : T_1| = 4$. Let
$\wt{\Sigma}_u$ be the group corresponding to  $\wt{\Sigma}$ in
$T_2$. Then  $|\wt{\Sigma}_u \cap T_1| \geq 2^6$. As any subgroup of
$\wt{\Sigma}$ of order at least $2^6$ is generated by conjugates of
$\wt{y}$, we have that $\wt \Sigma_u \cap T_1 \not \le \wt E$ by
Lemma~\ref{weakcl}. In particular, by the proof of Lemma~\ref{2conj}, $J>J_0$.  Therefore, we may suppose that there is a $\wt w \in \wt
\Sigma_u \cap T_1$ such that $\wt w$ inverts $\sigma$. Notice that
$(\wt \Sigma_u \cap T_1)\wt \Sigma$ is normal in $T_1\wt  \Sigma \in
\syl_2(\wt K_1)$. In particular, if $|(\wt \Sigma_u \cap T_1)\wt
\Sigma/\wt \Sigma|=2^2$, then $\wt w \wt \Sigma$ is centralized by a
maximal subgroup of $T_1\wt  \Sigma /\wt \Sigma$, which is
impossible. Hence $$|(\wt \Sigma_u \cap T_1)\wt \Sigma/\wt
\Sigma|\ge 2^3.$$ In particular, we have $|(\wt \Sigma_u \cap T_1
\cap \wt E)\wt \Sigma/\wt \Sigma|\ge 2^2$ and by Lemma~\ref{weakcl} all the non-trivial
elements are unitary transvections. This, however, contradicts
Proposition~\ref{fact} (viii) and proves the lemma.
\end{proof}

\begin{lemma}\label{weakcl2} We have $\wt{y}^{\wt{K}} \cap \wt K_1 \subseteq \wt{\Sigma}$. In particular, $\wt \Sigma$ is strongly closed in $\wt E$.
\end{lemma}

\begin{proof} Suppose that $\wt u \in \wt{y}^{\wt{K}} \cap \wt K_1 \setminus \wt{\Sigma}$. Then by Lemmas~\ref{weakcl} and \ref{2conj}, we get that $J > J_0$ and
  $\wt{u}$ inverts $\sigma$.
Furthermore,  all involutions in $\wt{\Sigma}\wt{u}$ are
conjugate. Hence, for $T_1\in \syl_2(C_{\wt{K}_1}(\wt{u}))$, we
have using Lemma~\ref{Normalizer} $|T_1| = 2^9$. Let $T_2$ be a Sylow 2-subgroup of
$C_{\wt{K}}(\wt{u})$ with $T_1 \leq T_2$ and $\wt{\Sigma}_u \leq
T_2$ be a $\wt K$-conjugate of $\wt{\Sigma}$ in $T_2$. \blue{ By Lemma \ref{tildex}, $(\wt{\Sigma}_u \cap T_1) \setminus \wt{\Sigma}$ does not contain elements $v$ with $|[v,\wt{\Sigma}]| = 4$. So by Lemma \ref{2conj} we have that $(\wt{\Sigma}_u \cap T_1)\wt{\Sigma}/\wt{\Sigma}$ inverts $\sigma\wt{\Sigma}$ and this yields
$\wt{\Sigma}_u \cap T_1 \subseteq \langle \wt u \rangle\wt \Sigma$.} Since $|\wt
\Sigma_u \cap T_1| \ge 2^5$, we now have that $\wt \Sigma_u \cap
T_1= \langle \wt u \rangle C_{\wt \Sigma}(\wt u)$ has order $2^5$. Hence
$T_2 = T_1\wt{\Sigma}_u$ and $T_2/\wt{\Sigma}_u \cong T_1/\langle\wt  u
\rangle C_{\wt \Sigma}(\wt u) \cong 2\times \Dih(8)$. But
$T_2/\wt \Sigma_u\cong \SDih(16)$  by Proposition \ref{fact} (ii) and we
thus have a contradiction. Hence $\wt{y}^{\wt{K}} \cap \wt K_1
\subseteq \wt{\Sigma}$.
\end{proof}

\begin{lemma}\label{strongclosed1} We have that $\wt{\Sigma}$ is strongly closed in $\wt{K}_1$.
\end{lemma}

\begin{proof} Assume by way of contradiction that there is some involution $\wt{u} \in \wt{K}_1 \setminus \wt{\Sigma}$, which is conjugate in $\wt{K}$ to
some element in $\wt{\Sigma}$.  By Lemma \ref{weakcl2} we have
$\wt{u} \sim_{\wt K} \wt{x}$. By Lemmas \ref{tildex} and
\ref{2conj} we have that $\tau \not= 1$ and we may assume that
$\wt{u}$ inverts $\sigma$. Furthermore we have $$C_{\wt
E/\wt{\Sigma}}(\wt u) \cong 2 \times \Sym(4).$$

Let $T_1$ be a Sylow 2-subgroup of $C_{\wt{K}_1}(\wt{u})$ and
$T_2$ be a Sylow 2-subgroup of $C_{\wt{K}}(\wt{u})$, which
contains $T_1$. Further let $\wt{\Sigma}_u$ be the normal subgroup
of $T_2$ which is $\wt K$-conjugate to  $\wt{\Sigma}$. Since, by
Proposition~\ref{fact} (ix), $C_{\wt \Sigma}(\wt u)$ is generated
by conjugates of $\wt y$, we have $C_{\wt{\Sigma}}(\wt{u}) \leq
\wt{\Sigma}_u$ by Lemma~\ref{weakcl2}.
 Since $(\wt \Sigma_u\cap T_1)\wt{\Sigma}/\wt{\Sigma} = \langle \wt
 u\rangle \wt{\Sigma}/\wt{\Sigma}$, we get
  $$T_3 = \wt{\Sigma}_u \cap T_1 = C_{\wt{\Sigma}}(\wt{u})\langle \wt{u} \rangle.$$
Therefore $T_3$ is normalized but not centralized by $\wt \Sigma$ and is centralized by
$\wt \Sigma_u$. \blue{We have that $\wt \Sigma$ and $\wt \Sigma_u$ are contained in $N_{\wt{K}}(T_3)$. Let $S_{\wt\Sigma}$ and $S_{\wt{\Sigma_u}}$ be Sylow 2-subgroups of $N_{\wt{K}}(T_3)$, which contain $\wt \Sigma$, $\wt{\Sigma}_u$, respectively.  As by
Lemma~\ref{K1Syl2} are  $\wt \Sigma$ is weakly closed in $S_{\wt\Sigma}$ and $\wt{\Sigma}_u$ is weakly closed in $S_{\wt{\Sigma_u}}$, we see that $\wt \Sigma$ and $\wt \Sigma_u$
are conjugate in $N_{\wt{K}}(T_3)$. But this is impossible as one centralizes $T_3$ and the other does not.}
\end{proof}

\begin{theorem}\label{structureK} We have $K = K_1$.
\end{theorem}

\begin{proof} Let $T \in \Syl_2(K)$. By  Lemmas~\ref{K1Syl2} and  \ref{strongclosed1} we have that $\wt{\Sigma}$ is strongly closed in $\wt T$ with respect to  $\wt{K}$.
Hence  an application of \cite{Goldschmidt} yields that $\wt{L} =
\langle \wt{\Sigma}^{\wt{K}} \rangle$ is an extension of a group of
odd order by a product of a 2-group and a number of  Bender groups.
Furthermore $\wt{\Sigma}$ is  the set of involutions in some Sylow
2-subgroup of $T\cap \wt{L}$. \blue{By Lemma \ref{5.4} we have that $K_1$ contains a Sylow 3-subgroup of $K$. As $C_K(\wt{\Sigma}) = \wt{\Sigma}$, we get that $O_{2^\prime}(\wt{L}) = O_{3^\prime}(\wt{L})$ now as $J_1$ normalizes $O_{2^\prime}(\wt{L})$ we get with Lemma \ref{Jsig}  that
$O_{2^\prime}(\wt{L}) = 1$.} As $\wt{K}_1$ acts primitively on $\wt{\Sigma}$,
either $\wt L = \wt{\Sigma}$ and we are done, or $\wt L$ is a  simple
group. So suppose that $\wt L$ is a simple group. Then  $N_{\wt
L}(\wt{\Sigma})$ acts transitively on $\wt{\Sigma}$, which is not
possible as $\Sigma$ is extraspecial. This proves that $K=K_1$.
\end{proof}

\section{Proof of the Theorem~\ref{MainTheorem}}

We continue with all the notation established in previous sections.
If $N_M(S)/S \cong \Dih(8)$, Theorem~\ref{MainTheorem}  follows with
Theorem \ref{NMSagain}. So we may assume that $N_M(S)/S \cong 2
\times \Dih(8)$. Using Theorem \ref{structureK} and
Lemma~\ref{containinK}  we get that $K/\Sigma \cong \Aut(\SU_4(2))$, $(3
\times \SU_4(2)){:}2$ or $\Sp_6(2)$.

Suppose that $K/\Sigma\cong \Sp_6(2)$.  Then \cite{Smith} implies
that  $G \cong \Co_2$ and consequently  $M=N_G(Z)$  has order
$2^8\cdot 3^6\cdot 5$ and  shape $3^{1+4}_+.2^{1+4}_-.\Sym(5)$, which is not
similar to a centralizer of type $\PSU_6(2)$ or $\F_4(2)$. This contradicts our
initial hypothesis.
 So suppose $K/\Sigma \cong \Aut(\SU_4(2))$ or  $(3 \times \SU_4(2)){:}2$.   Then Lemma \ref{notsimple} shows that $G$ possesses a subgroup $G_0$ of index two.
In particular we get $C_{G_0}(r)/\Sigma \cong \SU_4(2)$ or $3 \times
\SU_4(2)$. Now we see that $N_{G_0 \cap M}(S)/S \cong \Dih(8)$.
Hence Theorem \ref{NMSagain} gives $G_0 \cong \PSU_6(2)$ or
$\PSU_6(2){:}3$ and so  Theorem~\ref{MainTheorem} is proved.


\begin{thebibliography}{99}
\bibitem{Aschbacher} M.~Aschbacher, Finite Group Theory,
 Cambridge Studies in Advanced Mathematics, 10,
Cambridge University Press, Cambridge-New York, 1986.

\bibitem{BrauerSuzuki} Richard Brauer and    Michio Suzuki, On finite groups of
even order whose 2-Sylow group is a quaternion
 group,
 Proc. Nat. Acad. Sci. U.S.A.  45  1959 1757--1759.\bibitem{Atlas} J. H. Conway, R. T. Curtis, S. P. Norton, R. A. Parker and
R. A. Wilson, Atlas of finite groups. Maximal subgroups and ordinary
characters for simple groups,  Oxford University Press, 1985.
\bibitem{Goldschmidt} David Goldschmidt,  $2$-fusion in finite groups. Ann. of Math, 99 (1974), 70--117.
\bibitem{FeitThompson} Walter Feit and John G.  Thompson,   Finite groups which contain
a self-centralizing subgroup of order 3, Nagoya Math. J.  21  1962
185--197.
\bibitem{Gorenstein}  Daniel Gorenstein,  Finite groups,
Harper \& Row, Publishers, New York-London,   1968.
\bibitem{GLS2}Daniel  Gorenstein, Richard   Lyons and  Ronald
Solomon, The classification of the finite simple groups. Number 2,
Mathematical Surveys and Monographs, 40.2. American Mathematical
Society, Providence, RI,  1996.
\bibitem{GLS3}Daniel  Gorenstein, Richard   Lyons and  Ronald
Solomon, The classification of the finite simple groups, Number 3,
Mathematical Surveys and Monographs, 40.3. American Mathematical
Society, Providence, RI,  1998.

\bibitem{Hayden} John L. Hayden,   A characterization of the finite simple group
${\rm PSp}\sb{4}(3)$,
 Canad. J. Math.  25  (1973), 539--553.
\bibitem{Huppert} Bertram Huppert,   Endliche Gruppen,
Springer-Verlag, Berlin-New York  1967.

\bibitem{Isaacs}I. Martin. Isaacs,  Character theory of finite groups,  AMS Chelsea
Publishing, Providence, RI,  2006.
\bibitem{MSS} Ulrich  Meierfrankenfeld, Bernd  Stellmacher and  Gernot  Stroth,  The structure
theorem for finite groups with a large $p$-subgroup, preprint 2011.
\bibitem{Parker1} Chris
Parker,  A 3-local characterization of $\rm U\sb 6(2)$ and $\rm
Fi\sb
 {22}$,
 J. Algebra  300  (2006),  no. 2, 707--728.
\bibitem{McL}   Chris Parker  and Peter Rowley, A 3-local identification of the alternating group of degree eight,
the McLaughlin simple group and their automorphism groups,  Algebra 319, 2008, 1752 - 1775.
\bibitem{ParkerRowley2} Chris Parker and Peter  Rowley,  A $3$-local characterization of $\Co_2$,  J. Algebra  323  (2010),  no. 3, 601--621.

\bibitem{PS} Chris Parker and Gernot Stroth, An identification theorem for the sporadic simple groups $\F_2$ and $\M(23)$, preprint 2011.
\bibitem{PSStrong} Chris Parker and  Gernot  Stroth, Strongly $p$-embedded subgroups, Pure and Applied Mathematics Quarterly
Volume 7, Number 3
(Special Issue: In honor of
Jacques Tits)
797--858, 2011.
    \bibitem{PS1} Chris Parker and Gernot Stroth, Groups which are almost groups of Lie type in characteristic $p$, preprint 2011.
\bibitem{F42} Chris Parker and  Gernot  Stroth, ${\mathbf \F_4(2)}$ and its automorphism group, preprint 2011.
\bibitem{PSS} Chris Parker, M. Reza Salarian and Gernot Stroth, A characterisation of ${}^2\E_6(2)$, $\M(22)$ and $\Aut(\M(22))$ from a characteristic $3$ perspective, preprint 2011.
\bibitem{prince1} A. R. Prince,  A characterization of the simple groups ${\rm
PSp}(4,3)$ and ${\rm PSp}(6,2)$,  J. Algebra 45 (1977), no. 2,
306--320.
\bibitem{Seidel} Andreas Seidel,
Gruppen lokaler Charakteristik - eine Kennzeichnung von Gruppen vom Lie Typ in ungerader Charakteristik,
Dissertation, University Halle-Wittenberg, 2009.
http://digital.bibliothek.uni-halle.de/hs/content/titleinfo/397705.
\bibitem{Smith} Fredrick Smith,  A characterization of the $.2$ Conway Simple Group, J. Algebra 31, (1974), 91 -116.
\bibitem{Yoshida} Tomoyuki Yoshida,  Character-theoretic transfer, J. Algebra 52, (1978), 1-38.

\end{thebibliography}
\end{document}